\documentclass[11pt]{article}

\title{Oded Schramm: From Circle Packing to SLE}
\author{Steffen Rohde\footnote{University of Washington, Supported in part by NSF Grant DMS-0800968.}}
\date{July 12, 2010}

\input epsf.sty

\usepackage{amsfonts,amsthm,amsmath,latexsym,amssymb,graphicx}
\usepackage{subfigure}

\newtheorem{thm}{Theorem}[section]

\newif\ifdraft
\drafttrue

\def\b{{\beta}}
\def\g{{\gamma}}
\def\z{{\zeta}}
\def\zh{{\hat\zeta}}
\def\O{{\Omega}}
\def\D{{\mathbb D}}
\def\Z{{\mathbb Z}}
\def\R{{\mathbb R}}
\def\H{{\mathbb H}}
\def\C{{\mathbb C}}
\def\E{{\mathbb E}}
\def\P{{\mathbb P}}

\def\eps{\varepsilon}

\def\1{{\bf 1}}

\renewcommand{\Im}{\operatorname{Im}}

\begin{document}

\maketitle

\tableofcontents

\section{Introduction}

When I first met Oded Schramm in January 1991 at the University of California, San Diego, he 
introduced himself as a ``Circle Packer''. This modest description referred to his Ph.D. thesis around the  Koebe-Andreev-Thurston theorem and a discrete version of the Riemann mapping theorem, explained below. In a series of highly original papers, some joint with Zhen-Xu He, he created powerful new tools out of thin air, and provided the field with  elegant new ideas.
At the time of his deadly accident on September 1st, 2008, he was widely considered
as one of the most innovative and influential probabilists of his time. Undoubtedly, he is best known
for his invention of what is now called the Schramm-Loewner Evolution (SLE), and for his subsequent collaboration with Greg Lawler and Wendelin Werner that led to such celebrated results as a determination of the intersection exponents of two-dimensional Brownian motion and a proof of Mandelbrot's conjecture about the Hausdorff dimension of the Brownian frontier. But already his previous work bears witness to the brilliance of his mind, and many of his early papers contain both deep and beautifully simple ideas that deserve better knowing.

In this note, I will describe some highlights of his work in circle packings and the Koebe conjecture,
as well as on SLE. As Oded has co-authored close to 20 papers related to circle packings and more than 20 papers involving SLE, only a fraction can be discussed in detail here. 
The transition from circle packing to SLE was through a long sequence of influential papers
concerning probability on graphs, many of them written jointly with Itai Benjamini. 
I will present almost no work from that period (some of these results are described elsewhere in this volume, for instance in Christophe Garban's article on Noise Sensitivity). 
In that respect, the title of this note is perhaps misleading.

In order to avoid getting lost in technicalities,
arguments will  be sketched at best, and often ideas of proofs will be illustrated by analogies only.
In an attempt to present the evolution of Oded's mathematics, 
I will describe his work in essentially chronological order. 

Oded was a truly exceptional person: not only was his clear and innovative way of thinking 
an inspiration to everyone who knew him, but also his caring, modest and relaxed attitude 
generated a comfortable atmosphere. 
As inappropriate as it might be, I have included some personal anecdotes
as well as a few quotes from email exchanges with Oded, in order to at least hint at these
sides of Oded that are not visible in the published literature. 

This note is not meant to be an overview article about circle packings or  SLE.
My prime concern is to give a somewhat self-contained account of Oded's contributions.
Since  SLE has been featured in several excellent articles and even a book, 
but most of Oded's work on circle packing is accessible only through his original papers, 
the first part is a bit  more expository and contains more background. 
The expert in either field will find nothing new, 
and will find a very incomplete list of references.
My apology to everyone whose contribution is either unmentioned or, perhaps even worse, mentioned without proper reference. 

{\bf Acknowledgement:} I would like to thank Mario Bonk, Jose Fern\'andez, Jim Gill, Joan Lind, Don Marshall, Wendelin Werner and Michel Zinsmeister for helpful comments on a first draft. I would also like to thank Andrey Mishchenko for generating Figure \ref{fig2}, and Don Marshall for Figure \ref{incompfig}.

\section{Circle Packing and the Koebe Conjecture}
 Oded Schramm was able to create, seemingly without effort, ingenious new ideas and methods.
Indeed, he would be more likely to invent a new approach than to search the
literature for an existing one. In this way, in addition to proving
wonderful new theorems, he rediscovered many known 
results, often with completely new proofs. We will see many examples 
throughout this note.

\bigskip

Oded received his Ph.D. in 1990 under William Thurston's direction at Princeton.
His thesis, and the majority of his work until the mid 90's, was concerned with the 
fascinating topic of circle packings.
Let us begin with some background and a very brief overview of 
some highlights of this field 
prior to Oded's thesis. Other surveys are \cite{Sa} and \cite{Ste}.

\subsection{Background}

According to the Riemann mapping theorem, every {\it simply connected} planar domain, except the 
plane itself, is conformally
equivalent to a disc. The conformal map to the disc is unique, up to postcomposition
with an automorphism of the disc (which is a M\"obius transformation). The standard proof exhibits the map
as a solution of an extremal problem (among all maps of the domain {\it into} the disc, maximize the derivative at a given point). 
The situation is quite different for {\it multiply connected} domains, partly due to the lack of a standard target domain.
The standard proof can be modified to yield a conformal map onto a {\it parallel slit domain} (each complementary component is
a horizontal line segment or a point). Koebe showed that every 
{\it finitely connected} domain is conformally equivalent to a {\it circle domain}
(every boundary component is a circle or a point), 
in an essentially unique way. 
No proof similar to the standard proof of the Riemann mapping theorem is known. 

\begin{thm}[\cite{K1}]\label{koet} For every domain $\Omega\subset\C$ with finitely many
connected boundary components, there is a conformal map $f$ onto a 
domain $\Omega'\subset\C$ all of whose boundary components
are circles or points. Both $f$ and $\Omega'$ are unique up to a M\"obius transformation.
\end{thm}
\noindent
Koebe conjectured (p. 358 of \cite{K1}) that the same is true for infinitely connected domains. It later turned out that
{\it uniqueness} of the circle domain can fail (for instance, it fails whenever the set of point-components of the boundary has positive area, as a simple application of the measurable Riemann mapping theorem shows). But {\it existence} of a conformally equivalent circle domain is still open, and is known as Koebe's conjecture or ``Kreisnormierungsproblem''.
It motivated a lot of Oded's research.

\bigskip
There is a close connection between Koebe's theorem and circle packings.
A {\it circle packing} $P$ is a collection (finite or infinite) of closed discs $D$ in the 
two dimensional plane $\C$, or in the two dimensional sphere $S^2$, with disjoint interiors. 
Associated with a circle packing is its
{\it tangency graph} or {\it nerve} $G=(V,E)$,
whose vertices correspond to the discs, and such that two vertices are joined by an edge 
if and only if the corresponding discs are tangent. We will only consider packings whose tangency graph is connected.

\begin{figure}[htbp]
\centering
\includegraphics[width=120mm]{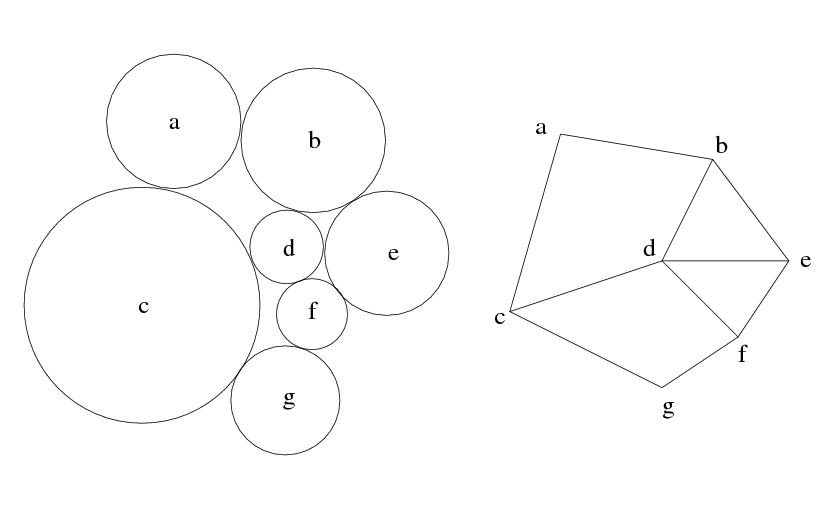}
\caption{\label{tanggraphfig} A circle packing and its tangency graph.}
\end{figure}

Conversely,  
the Koebe-Andreev-Thurston {\it Circle Packing Theorem} guarantees the existence of packings
with prescribed combinatorics. Loosely speaking, a planar graph is a graph that can be drawn in the plane so that edges do not cross. Our graphs will not have double edges (two edges with the same endpoints) or loops (an edge whose endpoints coincide). 
\begin{thm}[\cite{K2}, \cite{T}, \cite{A1}]\label{cpt} 
For every finite planar graph $G$, there is a circle packing in the plane with nerve $G$.
The packing is {\it unique} (up to M\"obius transformations) if $G$ is a triangulation of $S^2$.
\end{thm}
See the following sections for the history of this theorem, and sketches of proofs. In particular,
in Section \ref{ep} we will indicate how the Circle Packing Theorem \ref{cpt} can be obtained from 
the Koebe Theorem \ref{koet}, and conversely that the Koebe theorem
can be deduced from the Circle Packing Theorem.
Every finite planar graph can be extended (by adding vertices and edges as in Figure \ref{fig2c}) to
a triangulation, hence packability of triangulations implies packability of finite planar graphs 
(there are many ways to extend a graph to a triangulation, and uniqueness of the packing is no longer true).  The situation is more complicated for infinite graphs. Oded wrote several papers dealing with this case.
\begin{figure}[!htbp]
\centering
\includegraphics[width=120mm]{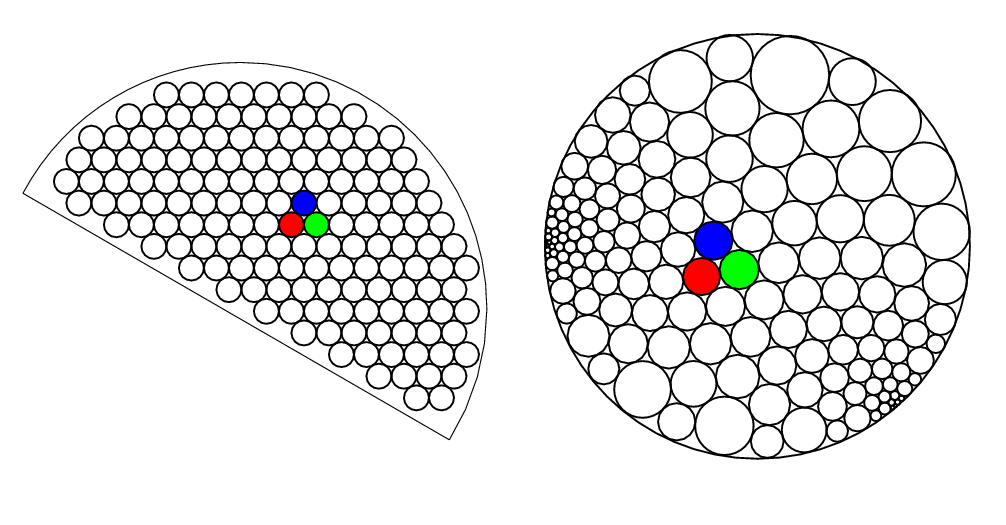}
\caption{\label{f1} A circle packing approximation to a Riemann map.}
\end{figure}
\bigskip
Thurston conjectured that circle packings approximate conformal maps, in the following sense:
Consider the {\it hexagonal packing} $H_{\eps}$ of circles of radius $\eps$ (a portion is visible in
Fig.~\ref{f1} and Fig.~\ref{fig2a}).
Let $\O\subset\C$ be a domain (a connected open set). 
Approximate $\O$ from the inside by a circle packing $P_{\eps}$ of circles of $\O\cap H_{\eps}$, as in
Fig.~\ref{f1} and Fig.~\ref{fig2a} (more precisely, take the connected component containing $p$ of the union of those circles whose six neighbors are still contained in $\O$). 
Complete the nerve of this packing by adding one vertex 
for each connected component of the complement to obtain a triangulation of the sphere (there are three new vertices $v_1, v_2, v_3$ in Fig.~\ref{fig2c}; the three copies of $v_3$ are to be identified). 
By the Circle Packing Theorem, there is a circle packing $P'_{\eps}$ of the sphere with the same tangency graph (Figures \ref{f1} and  \ref{fig2d} show these packings after stereographic projection from the sphere onto the plane; the circle corresponding to $v_3$ was chosen as the upper hemisphere and became the outside of the large circle after projection). 
Notice that each of the complementary components now corresponds to one (``large'') circle of $P'_{\eps}$, and the circles in the boundary of $P_{\eps}$ are tangent to these complementary circles.
Now consider the map $f_{\eps}$ that sends the centers of the circles of $P_{\eps}$ to the corresponding centers in $P'_{\eps}$, and extend it in a piecewise linear fashion.
\begin{figure}[!htbp]
\begin{center}
\subfigure[]{
	\label{fig2a}
	\includegraphics*[viewport=10 10 500 500, height=70mm]{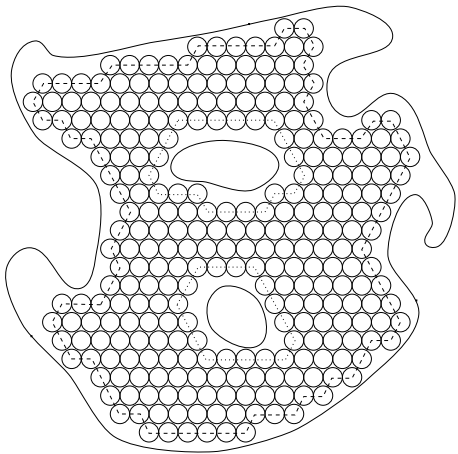}
}
\subfigure[]{
	\label{fig2b}
	\includegraphics*[viewport=10 10 500 500, height=70mm]{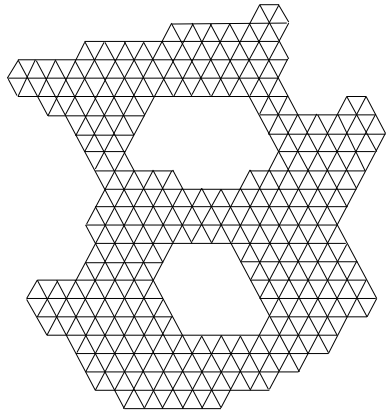}
}
\subfigure[]{
	\label{fig2c}
	\includegraphics*[height=90mm]{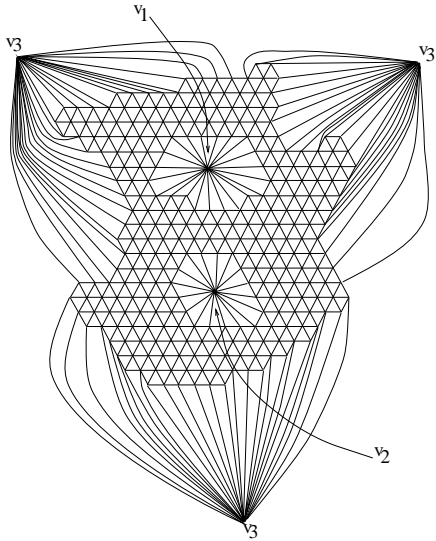}
}
\subfigure[]{
	\label{fig2d}
\includegraphics*[viewport=0 -30 250 250, height=80mm]{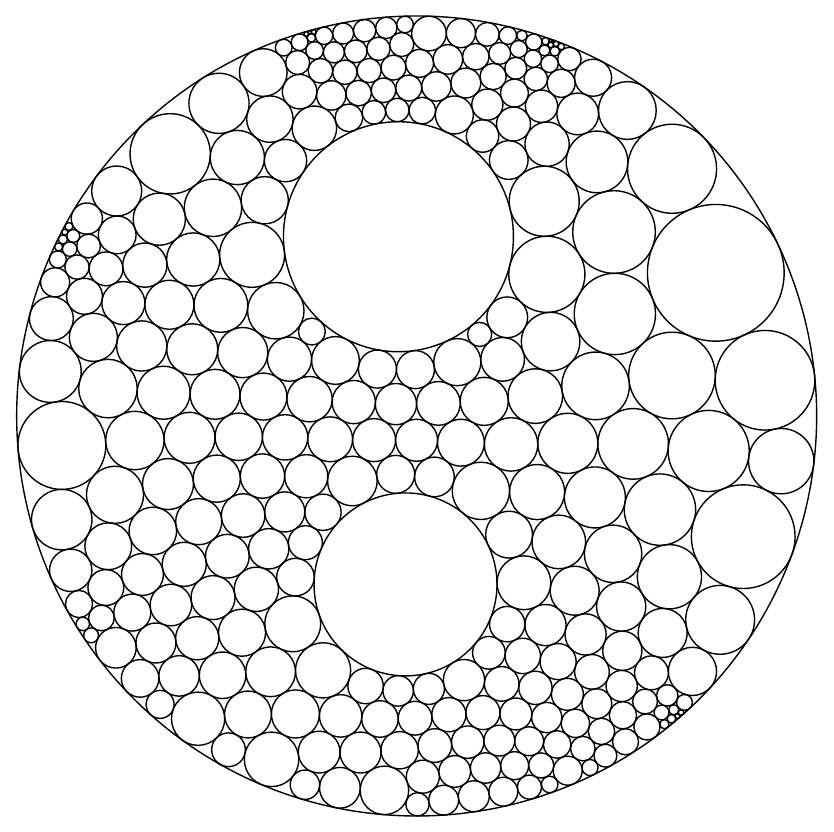}
}
\caption{A circle packing approximation of a triply connected domain,  its nerve, its completion to a triangulation of $S^2$, and a combinatorially equivalent circle packing; (a)-(c) are from Oded's thesis; thanks to Andrey Mishchenko for creating (d)}
\label{fig2}
\end{center}
\end{figure}
Rodin and Sullivan proved Thurston's conjecture 
that $f_{\eps}$ approximates the Riemann map,  if $\O$ is simply connected (see Fig.~\ref{f1}):
\begin{thm}\cite{RoSu}\label{rst} Let $\O$ be simply connected,  $p,q\in\O,$
and  $P'_{\eps}$  normalized such that the complementary circle is the unit circle, 
and such that the circle closest to $p$ (resp. $q$) corresponds to a circle containing $0$ 
(resp. some positive real number). 
Then 
the above maps $f_{\eps}$ converge to the conformal map $f:\O\to \D$ that is normalized by $f(p)=0$ and $f(q)>0,$ uniformly on compact subsets of $\O$ as $\eps\to0.$
\end{thm}
Their proof  depends crucially on the non-trivial {\it uniqueness} of the hexagonal packing  
as the only packing in the plane with nerve the triangular lattice.
Oded found remarkable improvements and generalizations of this theorem.
See Section \ref{convcp} for further discussion.

\subsection{Why are Circle Packings interesting?}

Despite their intrinsic beauty (see the book \cite{Ste} for stunning illustrations
and an elementary introduction), circle packings are interesting because they
provide a canonical and conformally natural way to embed a planar graph into a surface.
Thus they have applications to
combinatorics (for instance the  proof of Miller and Thurston \cite{MT}
of the Lipton-Tarjan separator theorem, see e.g.\ the slides of Oded's circle packing talk on his memorial webpage), 
to differential geometry (for instance the construction of minimal surfaces
by Bobenko, Hoffmann and Springborn \cite{BHS} and their references),
to geometric analysis (for instance,  the Bonk-Kleiner \cite{BK} quasisymmetric parametrization of
Ahlfors 2-regular LLC topological spheres)
to discrete probability theory (for instance,
through the work of Benjamini and Schramm on harmonic functions on graphs and recurrence on random planar graphs 
\cite{O12},\cite{O13}, \cite{O13b}) 
and of course to complex analysis (discrete analytic functions, conformal mapping). However,
Oded's work on circle packing did not follow any ``main-stream'' in conformal geometry or geometric function theory. 
I believe he continued to work on them just because he liked it. His interest never wavered,
and many of his numerous late contributions to Wikipedia were about this topic.

\bigskip
Existence and uniqueness are intimately connected. Nevertheless, for better readability I will discuss them in two separate sections.
\subsection{Existence of Packings}\label{ep}
\bigskip

 Oded applied the highest standards to his proofs  and 
was not satisfied with ``ugly'' proofs. 
As we shall see, he found four (!) different new
existence proofs for circle packings with prescribed combinatorics. 
Before discussing them, let us have a glance at previous proofs.

\bigskip

The Circle Packing Theorem was first proved by Koebe \cite{K2} in 1936. 
Koebe's proof of existence was based on his earlier result that every  planar domain $\O$ with
{\it finitely many} boundary components, say $m$, can be mapped  conformally onto a 
{\it circle domain}. A simple iterative algorithm, due to Koebe,
provides an infinite sequence $\O_n$ of domains conformally equivalent to $\O$ and such that $\O_n$ converges to a circle domain.
To obtain $\O_{n+1}$ from $\O_n$, just apply the Riemann mapping theorem
to the simply connected domain (in $\C\cup\{\infty\}$) containing $\O_n$ whose boundary 
corresponds to the $(n \mod m)-$th boundary component of $\O$.
With the conformal equivalence of  finitely connected domains and circle domains established,
a circle packing realizing a given tangency pattern can be obtained as a limit
of circle domains: Just construct a sequence of $m$-connected domains so that the boundary components approach each other  
according to the given tangency pattern. For instance, if the graph $G=(V,E)$ is embedded in the plane by specifying simple curves $\gamma_{e}:[0,1]\to S^2$, $e\in E,$ 
then the complement $\Omega_{\eps}$ of the set 
$$\bigcup_{e\in E} \gamma_e[0,1/2-\eps] \cup \bigcup_{e\in E} \gamma_e[1/2+\eps,1]$$
is such an approximation. 
It is not hard to show that the (suitably normalized) conformally equivalent circle
domains $\Omega'_{\eps}$  converge to the desired 
circle packing when $\eps\to 0.$

\bigskip
Koebe's theorem was nearly forgotten. In the late 1970's, Thurston rediscovered the circle packing theorem as an 
interpretation of a result of  Andreev  \cite{A1}, \cite{A2}
on convex polyhedra in hyperbolic space, and obtained uniqueness from
Mostow's rigidity theorem. He  suggested an algorithm to compute a circle packing 
(see \cite{RoSu}) and conjectured Theorem \ref{rst},
which started the field of circle packing. Convergence of Thurston's algorithm was proved in
\cite{dV1}.
Other existence proofs are based on a Perron family construction (see \cite{Ste}) and on a variational principle \cite{dV2}.

\bigskip
Oded's thesis \cite{O0} was chiefly concerned with a generalization of the existence theorem to packings with 
prescribed convex shapes instead of discs, and to applications. 
A consequence (\cite{O0}, Proposition 8.1) of his ``Monster packing theorem'' 
is, roughly speaking, that the circle packing theorem still holds if discs are replaced by smooth convex
sets.

\begin{thm}(\cite{O0}, Proposition 8.1)\label{convexcpt} 
For every  triangulation $G=(V,E)$ of the sphere, every $a\in V$,
every choice of smooth strictly convex sets $D_v$ for $v\in V\setminus\{a\},$
and every smooth simple closed curve $C$, there is a packing $P=\{P_v:v\in V\}$ with nerve
$G$, such that $P_a$ is the exterior of $C$ and each $P_v, v\in V\setminus\{a\}$ is
positively homothetic to $D_v$.
\end{thm}

%
%
%
\begin{figure}[!ht]
\centering
\includegraphics*[viewport=170 70 370 210, height=60mm]{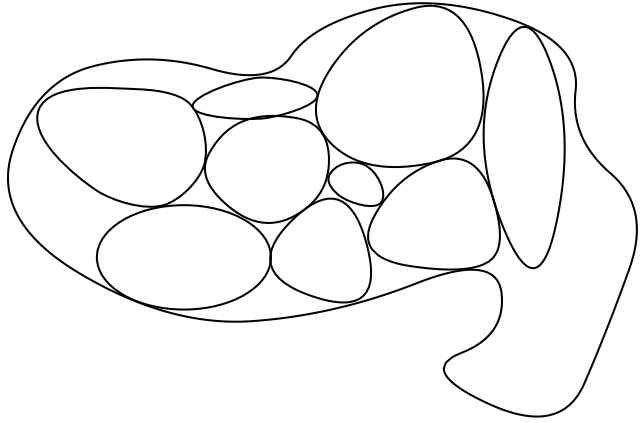}
\caption{\label{convexfig} A packing of convex shapes in a Jordan domain, from Oded's thesis}
\end{figure}

Sets $A$ and $B$ are positively homothetic if there is $r>0$ and $s\in\C$ with $A=r B+s.$
Strict convexity (instead of just convexity)
was only used to rule out that three of the prescribed sets could meet in one point
(after dilation and translation), and thus his packing theorem applied in much more generality.
Oded's approach was topological in nature: Based on a cleverly constructed spanning tree of $G$, 
he constructed what he called a ``monster''. This refers to a certain $|V|$-dimensional space of configurations of sets homothetic to the given convex shapes, with tangencies according to the tree, and certain non-intersection properties. 
Existence of a packing was then obtained as a consequence of Brower's fixed point theorem. Here is a poetic description, quoted from his thesis:

\bigskip
{\it One can just see the terrible monster swinging its arms in sheer rage, the tentacles causing a frightful hiss, as they rub against each other.}
\bigskip

Applying Theorem \ref{convexcpt} to the situation of Figure \ref{fig2}, with $D_v$ chosen as circles when $v\notin\{v_1, v_2, v_3\}$,
and arbitrary convex sets $D_{v_j}$, Oded adopted the Rodin-Sullivan convergence proof to obtain a new proof of 
the following generalization of Koebe's mapping theorem. The original proof of Courant, Manel and  Shiffman \cite{CMS} employed 
a very different (variational) method.

\begin{thm}(\cite{O0}, Theorem 9.1; \cite{CMS})\label{convexmapt} For every $n+1$-connected domain $\Omega$, every simply connected domain
$D\subset\C$ and every choice of $n$ convex sets $D_j$, there are sets $D_j'$ 
which are positively homothetic to $D_j$ such that
$\O$ is conformally equivalent to $D\setminus\cup_1^{n}D_j'$.
\end{thm}
Later \cite{O10} he was able to dispose of the convexity assumption, and proved the packing theorem for smoothly bounded but otherwise arbitrary shapes. As a consequence, he was able to generalize 
Theorem \ref{convexmapt} to arbitrary (not neccessarily convex) compact connected
sets $D_j$, thus rediscovering a theorem due to
Brandt \cite{Br} and Harrington \cite{Ha}. 

\bigskip Oded then developed a differentiable approach to the circle packing theorem.
In  \cite{O2b} he shows 
\begin{thm}(\cite{O2b},Theorem 1.1)\label{egg} Let $P$ be a 3-dimensional convex polyhedron, and let $K\subset \R^3$ be a smooth strictly convex body. Then there exists a convex polyhedron $Q\subset \R^3$ combinatorially
equivalent to $P$ which midscribes  $K.$
\end{thm}
Here ``$Q$ midscribes $K$'' means that all edges of $Q$ are tangent to $\partial K.$
He also shows that the space of such $Q$ is a six-dimensional smooth manifold, if the boundary of $K$ is smooth
and has positive Gaussian curvature. 
For $K=S^2,$ Theorem \ref{egg} has been stated by Koebe \cite{K2} and proved by Thurston \cite{T} using Andreev's
theorem \cite{A1}, \cite{A2}. Oded notes 
that  Thurston's midscribability proof based on the circle packing theorem can be reversed, so that Theorem \ref{egg} yields a new proof of the Circle Packing Theorem
(given a triangulation, just take $K=S^2$, $Q$ the midscribing convex polyhedron with the combinatorics of the packing, and for each vertex $v\in V$, let $D_v$ be the set of points on $S^2$ that are visible from $v$).

\bigskip
One defect of the continuity method in his thesis was that it did not provide a proof of uniqueness (see next section).
In \cite{O2} he presented a completely different approach to prove a far more general packing theorem, that had the added benefit of yielding uniqueness, too. A quote from \cite{O2}:

\bigskip
{\it It is just about the most general packing theorem of this kind that one could hope for (it is more general than I 
have ever hoped for).}

\bigskip
\noindent
A consequence of \cite{O2} (Theorem 3.2 and Theorem 3.5) is
%
%
\begin{thm} \label{ccpt}
Let $G$ be a planar graph, and for each vertex $v\in V$,  let ${\cal F}_v$ be a proper 3-manifold
of smooth topological disks in $S^2$, with the property that the pattern of intersection
of any two sets in ${\cal F}_v$ is topologically the pattern of intersection of two circles. Then there is a
packing $P$  whose nerve is $G$ and which satisfies $P_v\in {\cal F}_v$ for $v\in V$.
\end{thm}
\noindent
The requirement that ${\cal F}_v$ is a 3-manifold requires specification of a topology on the space of subsets of $S^2:$ Say that subsets $A_n\subset S^2$ converge to $A$  if
$\limsup A_n=\liminf A_n=A$ and $A^{c}=\rm{int}(\limsup A_n^c).$
An example  is obtained by taking a smooth strictly convex set $K$
in ${\bf R}^3$ and letting  ${\cal F}$ be the family of intersections $H\cap \partial K,$
where $H$ is any (affine) half-space intersecting the interior of $K.$ Specializing to $K=S^2$,
$\cal F$ is the familiy of circles and the choice ${\cal F}_v={\cal F}$ for all $v$ reduces to the circle packing theorem.

The proof of Theorem \ref{ccpt} is based on his 
{\it incompatibility theorem}, described in the next section.
It provides uniqueness of the packing (given some normalization), which is key to proving existence,
using continuity and topology (in particular invariance of domains).

\subsection{Uniqueness of Packings}\label{up}

\bigskip
I was always impressed by the flexibility of Oded's mind, in 
particular his ability to let go of a promising idea.
If an idea did not yield a desired result, it did not take long for him to come up with a 
completely different, and in many cases more beautiful, approach. 
He once told me that if he did not make
progress within three days of thinking about a problem,
he would move on to different problems. 

\bigskip Following Koebe and Schottky, 
{\bf uniqueness} of finitely connected circle domains (up to M\"obius images) is not hard to show, using the reflection principle: 
If two circle domains are conformally equivalent, the conformal map can be extended by reflection across 
each of the boundary circles, to obtain a conformal map between larger domains (that are still circle domains). 
Continuing in this fashion, one obtains a conformal map between complements of limit sets of reflection groups. 
As they are Cantor sets of area zero, the map extends to a conformal map of the whole sphere, hence is a M\"obius 
transformation. Uniqueness of the (finite) circle packing can be proved in a similar fashion. 
To date, the strongest rigidity result whose proof is based on this method is the following theorem of He and Schramm.
See \cite{B} for the related rigidity of Sierpinski carpets. 
\begin{thm}[\cite{O5}, Theorem A]\label{to5} If $\O$ is a circle domain whose boundary has $\sigma-$finite length, then $\O$ is rigid
(any conformal map to another circle domain is M\"obius).
\end{thm}

\bigskip
For {\it finite} packings, there are several technically simpler proofs. 
The shortest and most elementary of them is deferred to the end of this section, since I believe it has been discovered last. Rigidity of {\it infinite} packings lies deeper. The rigidity of the hexagonal packing,
crucial in the proof of the Rodin-Sullivan theorem as elaborated in Section \ref{convcp} below, was originally
obtained from deep results of Sullivan's concerning hyperbolic geometry. He's thesis \cite{He} gave a quantitative 
and simpler proof, still using the above reflection group arguments and the theory of quasiconformal maps.
In one of his first papers \cite{O1}, Oded gave an elegant combinatorial proof that at the same time was more general:
\begin{thm}[\cite{O1}, Theorem 1.1]\label{cprigidity} Let $G$ be an infinite, planar triangulation and
$P$ a circle packing on the sphere $S^2$ with nerve $G.$ If $S^2\setminus {\rm carrier}(P)$ is at most countable, then $P$ is rigid (any other circle packing with the same combinatorics is M\"obius equivalent).
\end{thm}
The {\it carrier} of a packing $\{D_v:v\in V(G)\}$
is the union of the (closed) discs $D_v$ and the ``interstices'' (bounded by three mutually touching circles) in the complement of the packing. 
The rigidity of the hexagonal packing 
follows immediately, since its carrier is the whole plane.

The ingenious new tool is his {\it Incompatibility Theorem}, a combinatorial analog to the conformal modulus
of a quadrilateral. To fully appreciate it, lets first look at its classical continuous counterpart,
and defer the statement of the Theorem to Section \ref{sincomp} below.

\subsubsection{Extremal length and the conformal modulus of a quadrilateral}\label{modulus}

If you conformally map a 3x1-rectangle to a disc, such that the center maps to the center,
what fraction of the circle does the image of one of the two short sides occupy?
Despite having known the effect of ``crowding'' in numerical conformal mapping, I was surprised
to learn of the numerical value of $0.0114...$ from Don Marshall (see \cite{MS}.)
Of course, the precise value can be easily computed as an elliptic integral, but if asked for a
rough guess, most answers are around 1/10 (the uniform measure with respect to length would give 1/8). 
Oded's answer, after a moments thought (during a tennis match in the early 90's), was 1/64, 
reasoning that this is the probability of a planar random walker to take each of his first three steps 
``to the right''. 

\bigskip
An important classical conformal invariant, masterfully employed by Oded in many of his papers, is the 
{\it modulus} of a quadrilateral. Let $\O$ be a simply connected domain in the plane that is bounded by a simple closed curve, and let $p_1,p_2,p_3$ and $p_4$ be four consecutive points on $\partial\O.$
Then there is a unique $M>0$ such that there is a conformal map $f:\O\to[0,M]\times[0,1]$ and such that
$f$ takes the $p_j$ to the four corners with $f(p_1)=0$ (by a classical theorem of Caratheodory, $f$ extends
homeomorphically to the boundary of the domains). There are several quite different instructive
proofs of uniqueness of $M$. Each of the following three techniques has a counterpart in the circle
packing world that has been employed by Oded.
Suppose we are given two rectangles and a conformal map $f$ between them taking corners to corners.

\bigskip
One method to prove uniqueness is to repeatedly reflect $f$ across the sides of the rectangles. The resulting extention is a conformal map of the plane, hence linear, and it follows that the aspect ratio is unchanged. This is similar to the aforementioned Schottky group argument.

\bigskip
A second method is to explicitly define a quantity $\lambda$ depending on a configuration 
$(\O,p_1,...,p_4)$
in such a way that it is conformally invariant and such that one can compute 
$\lambda$ for the rectangle $[0,M]\times[0,1]$. This is achieved by the {\it extremal length} of the
family $\Gamma$ of all rectifiable curves $\gamma$ joining two opposite ``sides'' $[p1,p2]$ and
$[p3,p4]$ of $\O.$ The extremal length of a curve family $\Gamma$ is defined as 
\begin{equation}\label{el}
\lambda(\Gamma)=\sup_{\rho}\frac{(\inf_\gamma  \int_{\gamma} \rho |dz| )^2 }{  \int_{\C} \rho^2 dx dy},
\end{equation}
where the supremum is over all ``metrics'' (measurable functions) $\rho:\C\to[0,\infty)$.
For the family of curves joining the horizontal sides in the rectangle $[0,M]\times[0,1]$,
it is not hard to show $\lambda(\Gamma)=M.$
This simple idea is actually one of the most powerful tools of geometric function theory.
See e.g. \cite{P} or \cite{GM} for references, properties and applications.

\bigskip
Discrete versions of extremal length (or the ``conformal modulus'' $1/\lambda$) have
been around since the work of Duffin \cite{Duf}. In conformal geometry, they 
have been very succesfully employed beginning with the groundbreaking paper  \cite{C}.   
Cannon's extremal length on a graph $G=(V,E)$ 
is obtained from \eqref{el} by viewing non-negative functions $\rho:V\to[0,\infty)$
as metrics on $G$, defining the length of a ``curve'' $\gamma\subset V$ 
as the sum $\sum_{v\in\gamma} \rho(v)$, and the ``area'' of the graph as $\sum \rho(v)^2.$
See \cite{CFP1} for an account of Cannon's discrete Riemann mapping theorem,
and for instance the papers \cite{HK} and \cite{BK} concerning applications to quasiconformal geometry.
Oded's applications to square packings and transboundary extremal length are briefly discussed
in Section \ref{other} below.

\begin{figure}[!htbp]
\centering
\includegraphics*[viewport=10 200 580 620, width=70mm]{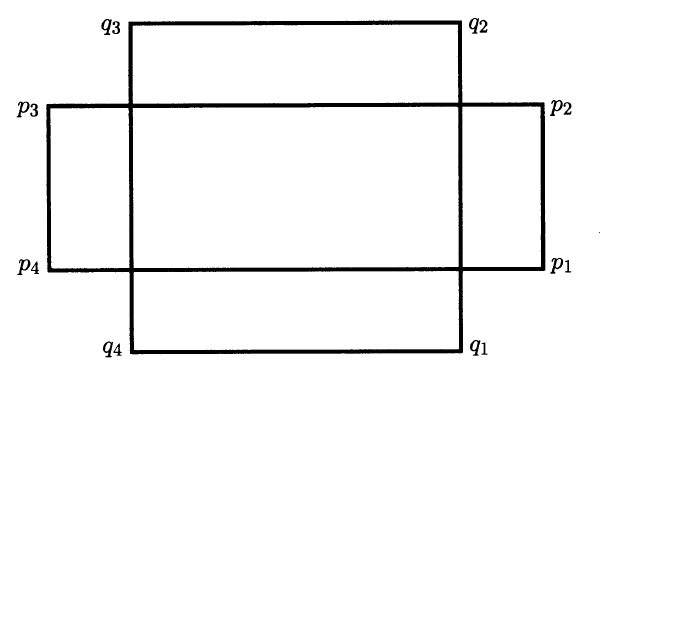}
\caption{\label{rect} Conformally inequivalent rectangles; from \cite{O3}.}
\end{figure}

\bigskip
A third and very different method is topological in nature and is one of the key ideas in \cite{O3}.
Suppose we are given two rectangles $\O,\O'$ with different aspect ratio and overlapping as in 
Fig.~\ref{rect}, 
and a conformal map $f$ between them mapping  corners to corners. Then the difference $f(z)-z$ is
$\neq0$ on the boundary $\partial\O$. Traversing $\partial\O$ in the positive direction, inspection of Fig.~\ref{rect} shows that the image curve under $f(z)-z$ winds around 0 in the negative direction.
But a negative winding is impossible for analytic functions (by the argument principle, the winding number counts the number of preimages of $0$).
\subsubsection{The Incompatibility Theorem}\label{sincomp}
\begin{figure}[htbp]
\centering
\includegraphics*[viewport= 0 0 1100 350, width=170mm ]{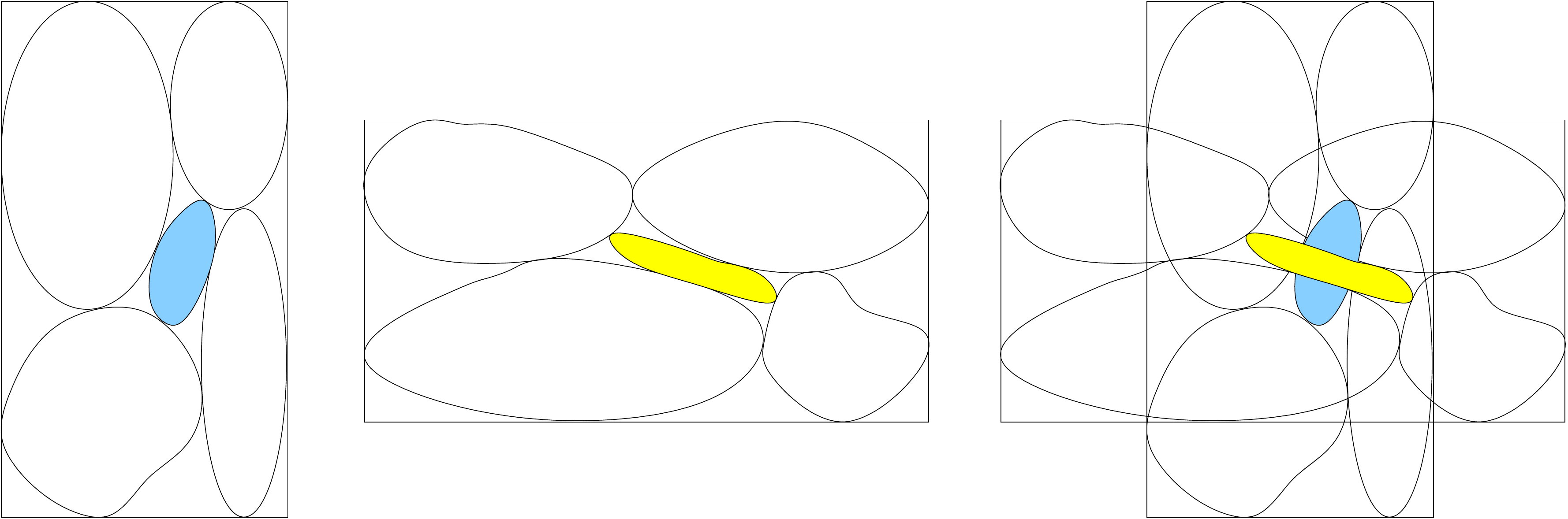}
\caption{\label{fincomp} An incompatibility at the center.}
\label{incompfig}
\end{figure}

Again consider the overlapping rectangles $\O,\O'$ of Fig.~\ref{rect}, 
and  two combinatorially equivalent
packings $P,P'$ whose nerves triangulate the rectangles, as in  Fig.~\ref{fincomp}. 
Assume for simplicity that the sets $D_v$ and $D'_v$ of the packings are closed
topological discs (except for the four sides $D_1,...D_4$, $D_1',...,D_4'$
of the rectangles, which are considered to be sets of the packing).
Intuitively, two topological discs $D$ and $D'$ are called {\it incompatible} if they intersect as in Fig.~\ref{rect}.
More formally, say that $D$ {\it cuts} $D'$ if there are two points in $D'\setminus {\rm interior}(D)$
that cannot be connected by a curve in ${\rm interior}(D'\setminus D)$. Then Oded calls $D$ and $D'$  incompatible if $D$ cuts $D'$ or $D'$ cuts $D.$ 
As he notes, {\it the motivation for the definition comes from the simple but very important observation that the possible patterns of intersection
of two circles are very special, topologically.} Indeed, any two circles are compatible.
\begin{thm}[\cite{O1}, Theorem 3.1]\label{incompatibility}
There is a vertex $v$ for which $D_v$ and $D'_v$ are incompatible.
\end{thm}

Oded calls this result a combinatorial version of the modulus.
However, it has rather little in common with the above notion of discrete modulus, except
for the setup. 

Oded's clever proof by induction on the number of sets in the packing uses arguments from plane topology.
An immediate consequence is that two rectangles cannot be packed by the same circle pattern,
unless they have the same modulus $M$ and hence are similar:
if they could, just place the two packings on top of each other as in Fig.~\ref{fincomp} and obtain 
two incompatible circles, a contradiction. In the same vein, it is not
difficult to reduce the proof of the rigidity 
Theorem \ref{cprigidity} to an application of the incompatibility theorem.

\subsubsection{A simple uniqueness proof}

\bigskip
To end this section, here is a beautifully simple proof of the rigidity of finite circle packings
whose nerve triangulates $S^2$. I copied it from the wikipedia (search for circle packing theorem), and  believe it is due to Oded. As before, stereographically project the packing
to obtain a packing of discs in the plane. This time, assume that the north pole belongs to the complement of the discs, so that the planar packing will consist of three ``outer'' circles and the remaining circles contained in the interstice between them.

\bigskip
{\it ``There is also a more elementary proof based on the maximum principle, which we now sketch. The key observation here is that if you look at the triangle formed by connecting the centers of three mutually tangent circles, then the angle formed at the center of one of the circles is monotone decreasing in its radius and monotone increasing in the two other radii. Consider two packings corresponding to G. First apply reflections and 
M\"obius transformations to make the outer circles in these two packings correspond to each other and have the same radii. Next, consider a vertex v where the ratio between the corresponding radius in the one packing and the corresponding radius in the other packing is maximized. Since the angle sum formed at the center of the corresponding circles is the same (360 degrees) in both packings, it follows from the above observation that the radius ratio is the same at all the neighbors of v as well. Since G is connected, we conclude the radii in the two packings are the same, which proves uniqueness.''}
\subsection{Koebe's Kreisnormierungsproblem}
Koebe's 1908 conjecture \cite{K1} that every planar domain can be mapped conformally onto a circle domain is still open, despite considerable effort by Koebe and others. Important contributions were made by Gr\"otzsch, Strebel, Sibner and others.
One difficulty is the aforementioned lack of uniqueness. Another problem is that Theorem \ref{convexmapt} is not true in the infinitely connected case, as the following
example from \cite{O9} illustrates: If 
$K=\{x+i y: x=0, \pm1, \pm\frac12, \pm\frac13,..., y\in[-1,1]\}$, 
and if $D=\hat\C\setminus K$,
then there is no conformal map $f$ of $D$, normalized by $f(z)-z\to 0$ as $z\to\infty,$
such that the component $\{i y: y\in[-1,1]\}$
of $\partial D$ corresponds to a horizontal line segment (or a point) while the other complementary components of $f(D)$ are vertical line segments. The same example also illustrates the fundamental continuity problem: There is a circle domain $D'$
conformally equivalent to D, but the boundary component corresponding to $\{i y: y\in[-1,1]\}$ is just a point, so that the conformal
map from $D'$ to $D$ cannot be extended to the boundary.

\bigskip
The first joint paper of  He and Schramm provided a breakthrough:
\begin{thm}[\cite{O3}]\label{countable}
If $\O$ has at most countably many boundary components, then $\O$ is conformally equivalent 
to a circle domain $\O'$, and $\O'$ is unique up to M\"obius transformation.
\end{thm}
Essentially, this result is still the strongest to date. 
Oded later \cite{O9} gave a conceptually different and  simpler proof based on 
his transboundary extremal length,  which also applies to certain classes of domains with uncountably many boundary components. 

\bigskip
The proof in \cite{O3} used transfinite induction and was based on the topological concept of the 
{\it fixed-point index}. I will illustrate the beautiful idea by sketching their proof of uniqueness.
As it turned out, this argument for uniqueness had been given earlier by Strebel \cite{Str}.
The simple but crucial idea is to use the following (see \cite{O3}, Lemma 2.2): If $f$ is a fixed-point free 
orientation preserving homeomorphism between two  circles $C'$ and $C''$, 
then the winding number of the curve
$f(z)-z, z\in C',$ around 0 is non-negative (recall Fig.~\ref{rect} for a situation where the winding
number is negative).
Let $f:\O'\to \O''$ be a conformal map and 
assume for simplicity  that $f$ extends continuously to the boundary
(in case of finitely many boundary components this is immediate from the reflection principle, but in
the countable case this step is non-trivial), and that $f$ has no fixed points on the boundary. Composing with M\"obius transformations, 
we may assume that $\infty\in\O'$ and that $f(z)=z+a_1/z+a_2/z^2+\cdots.$ We want to show that $f$ is
the identity. If not, denote $a_j$ the first
non-zero Taylor coefficient, then $f(z)-z$ has winding number  $-j$ as $z$ traverses a large circle $|z|=R,$ because $f(z)-z$ behaves like $a_j z^{-j}$.
Moreover, each circular boundary component maps to a circular component. These boundary components 
are oriented negatively (to keep the domain to the left) and thus, by the above crucial idea, contribute a non-positive number to the winding of $f(z)-z, z\in\partial(\O\cap \{|z|\leq R\}$ around 0.
Hence the total winding number is negative, contradicting their generalization of 
the argument principle (the winding number counts the number of zeroes of $f(z)-z.$) Of course, I have swept most details under the rug, most notably the proof of continuity based on a powerful generalization of Schwarz' Lemma to circle domains (Theorem 0.6 in  \cite{O3}).

\bigskip
Combining the fixed-point index method of \cite{O3} with an analysis of quasiconformal deformations
using the reflection group approach and Sullivan's rigidity theorems, He and Schramm \cite{O7}
improve Theorem \ref{countable} to domains $\O$ for which all boundary components are circles or points
except those in a countable and closed family. They also obtain the following
generalization of the Riemann mapping theorem. Let $A\subseteq\C$ be simply connected. 

\begin{thm}[\cite{O7},\cite{O8}]\label{relative} If $\O\subset A$ is a relative circle domain
(each connected component of $A\setminus\O$ is a  point or a closed disc), then there is a 
relative circle domain $\O^*$ in $\D$ conformally equivalent to $\O$, and so that $\partial A$
corresponds to $\partial\D.$
Conversely, if $\O^*$ is a relative circle domain in $\D$, there is such $\O\subset A.$
\end{thm}
The converse direction
is the main result of \cite{O8}.

\subsection{Convergence to conformal maps}\label{convcp}
Let us return to the setting of the Rodin-Sullivan Theorem \ref{rst} about convergence
of the discrete map $f_{\eps}$ to the conformal map $f.$
Consider the piecewiese linear extension of $f_{\eps}$ from the carrier of $P$ 
to the carrier of $P'$ that maps equilateral triangles to the corresponding triangles (formed
by the centers of $P'$). By the elementary ``Ring Lemma'' of \cite{RoSu}, the angles of these
triangles are bounded away from $0$ and $\pi$ (so that $f_{\eps}$ is quasiconformal with dilation
uniformly bounded above). 
At the heart of the Rodin-Sullivan proof is the {\it uniqueness} of the hexagonal packing  as the only 
packing in the plane with nerve the triangular lattice (see the discussion in Section \ref{up}).
It rather easily implies that tangent circles
centered in a compact set of $\O$ correspond to tangent circles in $\D$ whose radii are asymptotically equal as $\eps\to0.$ Hence the triangles in $P'$ are nearly equilaterals when $\eps$ is small
(the angles tend to $\pi/3$), so that $f_{\eps}$ is nearly angle preserving in each triangle.
Now the theory of quasiconformal maps readily yields equicontinuity of the family of maps 
$f_{\eps}$, and shows that every subsequential limit $\lim f_{\eps_j}$
is a conformal map. The theorem follows from uniqueness of normalized conformal maps. 

He's thesis \cite{He} provided a quantitative estimate for the rate of convergence of the angles
(the difference to $\pi/3$ is $O({\eps}$)). This estimate was known to imply
convergence of the ratio of corresponding radii $\rm{rad}(D')/\rm{rad}(D)$ to the absolute value
$|f'|$ of the derivative. A probabilistic proof of $C^0$ (locally uniform) convergence of circle packings was given by Stephenson \cite{Step}.
Convergence of $f_{\eps}$ to $f$ for packings other than the hexagonal
was proved in \cite{HR}, under the assumption of bounded valency of the graph.
In \cite{DHR}, the quality of convergence was improved to convergence in $C^2$ (that is,
convergence of first and second derivatives; strictly speaking,
instead of $f_{\eps}$ they considered the ``piecewise M\"obius'' map that sends interstices between triples of mutually tangent circles to the corresponding interstices). He and Schramm \cite{O11} found
an elementary new convergence proof, based on the topological ideas discussed above and thus avoiding 
quasiconformal maps. Their proof also gave convergence up to $C^2$, and worked in a more general 
setting. In particular, it does not need the assumption of uniformly bounded degree of \cite{HR}.

In the remarkable paper \cite{O16}, He and Schramm proved 
$C^\infty$-convergence of hexagonal disk packings to the Riemann map:

\begin{thm}[\cite{O16}, Theorem 1.1]
The discrete functions $f_{\eps}: V_{\eps}\to\D$ converge in $C^\infty$
to the Riemann mapping $f: \O\to\D$, in the sense that the 
discrete partial derivatives of $f_{\eps}$ of any order converge
locally uniformly to the corresponding partial derivatives of $f$.
\end{thm}
The discrete first-order derivatives for $v\in V_{\eps}$ are
$$\partial_{\eps,k} f_{\eps} (v) = \eps^{-1}(f_{\eps}( v + \eps\ \omega^k)-f_{\eps} (v)),$$
where $k\in{0,1,...,5}$ and $\omega=(1+i\sqrt 3)/2$ is a $6-$th root of unity.
In particular, it follows that $(\partial_{\eps,0})^k f_{\eps}$ converges to the $k-$th
derivative $f^{(k)}$ locally uniformly on $G$. 

The Schwarzian derivative
\begin{equation}\label{schwarzian}
S(f)(z) = \frac{f'''(z)}{f'(z)}-\frac32\frac{f''(z)^2}{f'(z)^2}
\end{equation}
of a locally univalent analytic function measures the deviation of $f$ from a M\"obius transformation,
in particular $S(f)\equiv0$ if $f$ is M\"obius.
A key idea in the proof is to define a discrete
analog of the Schwarzian derivative, to compute the (discrete) Laplacian of this Schwarzian,
and to employ a regularity theorem for discrete elliptic equations to obtain boundedness of
all partials of the Schwarzian. The definition of the discrete Schwarzian is the circle packing analog
of an invariant that Oded so masterfully employed in his earlier work  \cite{O14} 
on circle patterns with the combinatorics of the square grid.

\subsection{Other topics}\label{other}

Oded's approach to both mathematics and to life was extraordinarily innovative 
and unaccepting of conventions. Notions that most people take for granted
without even thinking about, he would open-mindedly question, often coming
up with amazing alternative solutions. 
For example, I would not even think about camping  on the foot of a glacier without
a sleeping bag. Climbing little Tahoma peak with Oded, he proved to me that
even this idea can be pursued. It was perhaps one of his less successful innovations, though.

\begin{figure}[!htbp]
\begin{center}
\subfigure[]{
	\label{figsqa}
	\includegraphics*[viewport=-100 -100 500 500, height=60mm]{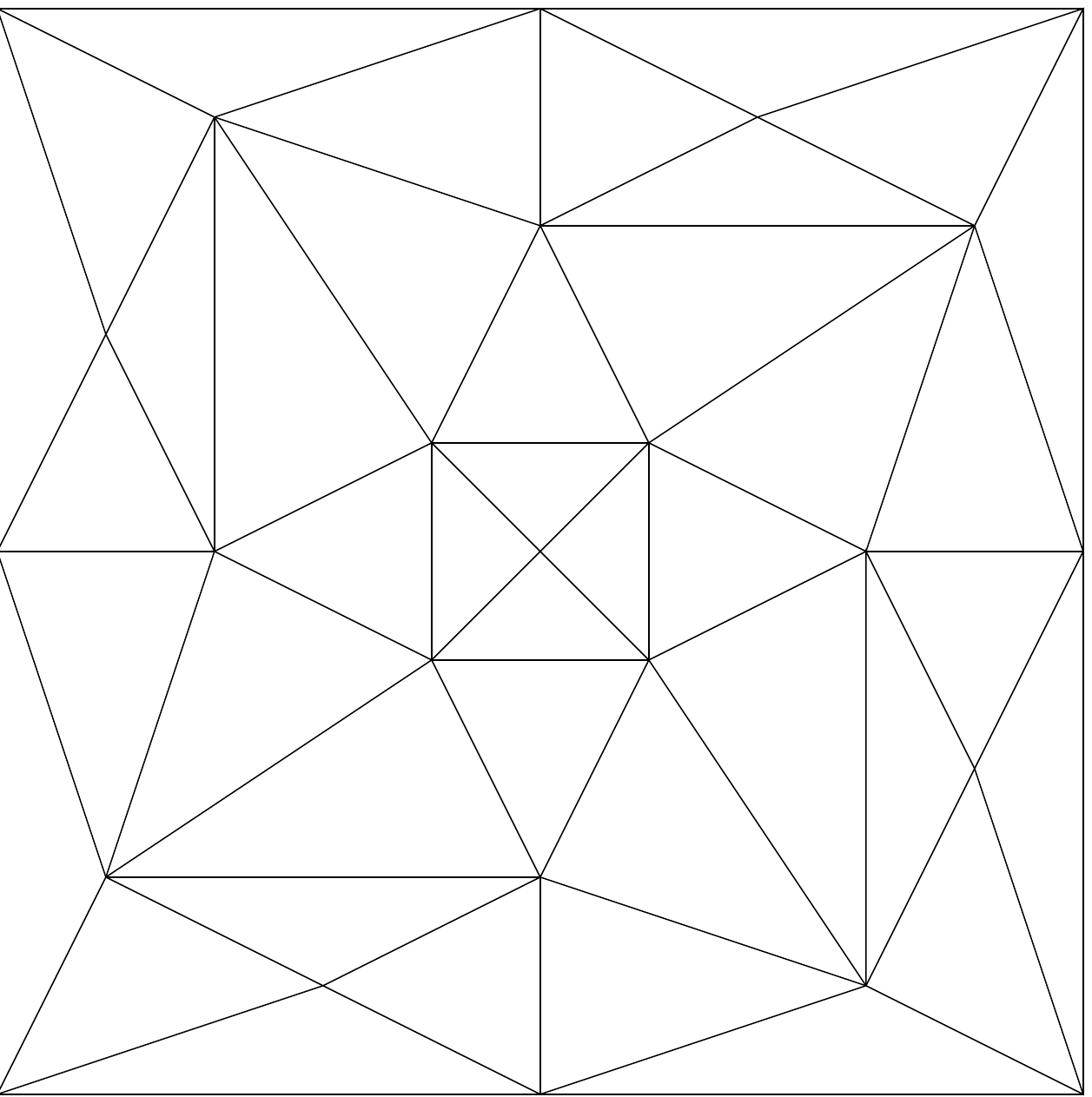}
}
\subfigure[]{
	\label{figsqb}
	\includegraphics*[viewport=10 10 600 600, height=60mm]{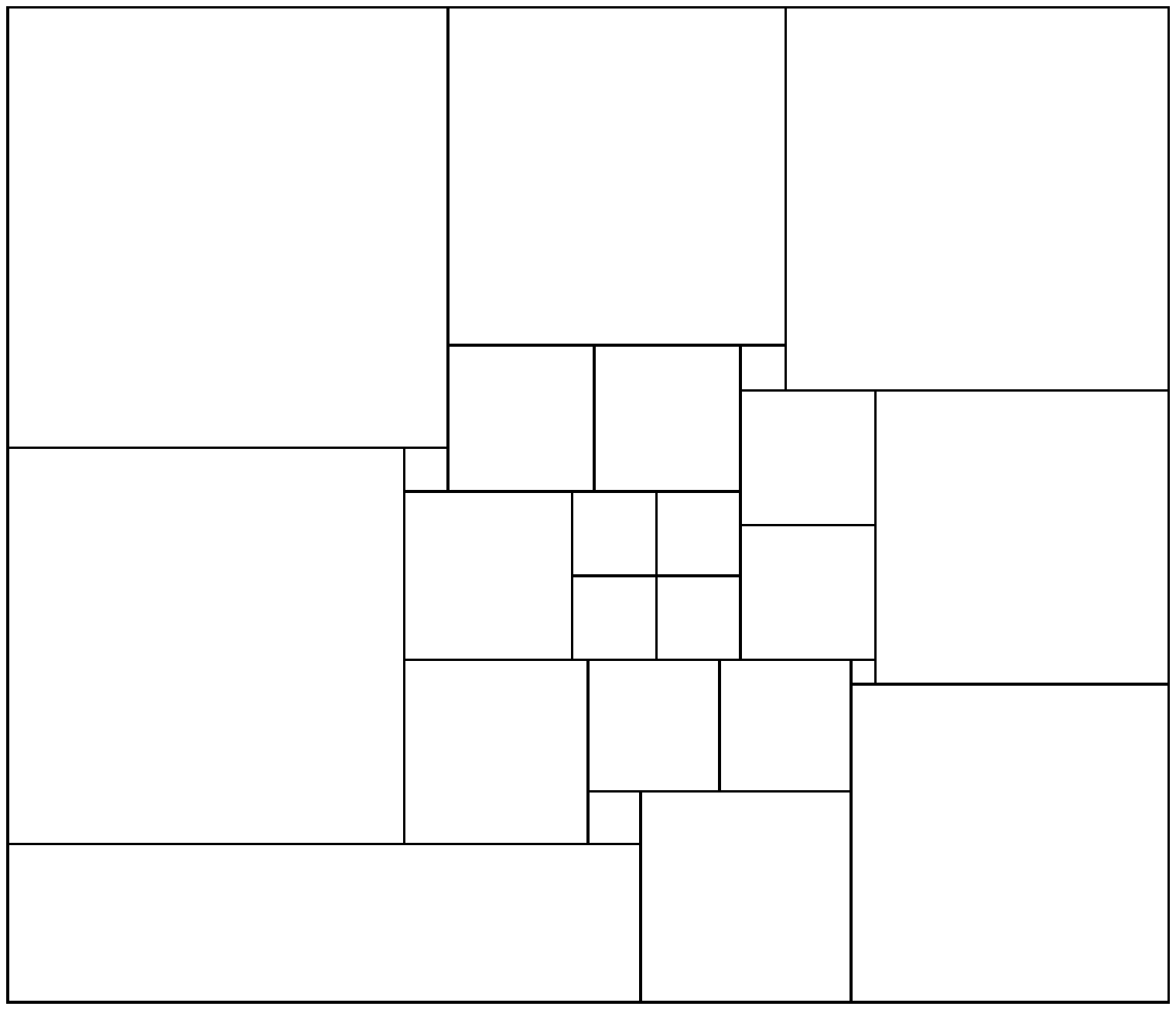}
}
\caption{A triangulation and the associated square packing. Thanks to David Wilson for providing this
figure from  \cite{O4}}
\label{figsq}
\end{center}
\end{figure}

\bigskip
In the lovely paper \cite{O4}, Oded shows that for each triangulation
$G$ of a quadrilateral, there is  a packing of a rectangle $R$ by (horizontal) squares with the combinatorics
of $G$ (a square might degenerate to a point, as in Figure \ref{figsq}). 

\noindent
The packing is actually a tiling: Indeed,
Oded points out the following simple observation.

\bigskip
\noindent
{\it Let $P_a$, $P_b,$ $P_c$ be three rectangles whose edges are  parallel to the coordinate axis. Suppose that the
intersection of every two of these rectangles is nonempty. Then $P_a\cap P_b\cap P_c \neq \emptyset.$
}

\bigskip
\noindent
The same tiling theorem was obtained independently by Cannon, Floyd and Parry in \cite{CFP1}. 
Both  employ Cannon's discrete extremal length (see Section \ref{modulus}) and obtain the side lengths $s(v)$
of the squares as the weights $\rho(v)$ of the extremal metric
(corresponding to the family of ``combinatorial curves'' joining two opposite sides of the quadrilateral).
It is quite different from the classical square packings of Brooks, Smith, Stone and Tutte \cite{BSST}, 
in particular, since the metrics considered here live on the vertices rather than the edges of the graph.

\bigskip
A very similar idea is exploited in the important paper \cite{O9}. 
The classical setting of extremal length
(recall \eqref{el}) is  a family $\Gamma$ of curves contained in a domain $\O.$ Invariance
$$\lambda (f(\Gamma)) = \lambda(\Gamma)$$ under conformal maps of $\O$ is almost trivial 
(just pull back metrics from $f(\O)$). Oded's notion of {\it transboundary extremal length}
$\lambda_{\Omega}(\Gamma)$ applies to curve families $\Gamma$ that are not necessarily contained in
$\O.$ The metrics are now replaced by generalized metrics  $\rho$ that, roughly speaking, also assign length to complementary components. The length 
$\int_{\gamma} \rho |dz|$ is replaced by 
$\int_{\gamma\cap \O} \rho |dz|+\sum_{p} \rho(p)$ if $\gamma$ is not contained in $\O,$
where the sum is over all boundary components of $\O$ that $\gamma$ meets. Then the definition is
$\lambda_{\Omega}(\Gamma)= \sup_{\rho} (\inf_\gamma \int_{\gamma\cap \O} \rho |dz|+\sum_{p} \rho(p))/
\int_{\C} \rho^2 dx dy$, and conformal invariance is again immediate.
Using this innocent looking extension, Oded provides an elegant self-contained proof of the countable Koebe conjecture, and moreover is able to deal with the case 
of domains for which the complementary component
satisfy a certain fatness condition (${\rm area}(A\cap B(x,r))\geq c r^2$ for each component $A$,
each $x\in A$ and each disc $B(x,r)$ that does not contain $A$).

\bigskip
Circle packings corresponding to infinite graphs $G$ can be obtained by taking Hausdorff limits 
of packings corresponding to finite subgraphs, but where do they ``live''? Beardon and Stephenson \cite{BSt1},\cite{BSt2}
have shown, under the assumption that the degrees of the vertices are uniformly bounded,
that the carrier of such a packing is either the plane (call this case {\it parabolic}),
or that it can be chosen to be the disc ({\it hyperbolic}). They also showed that both cases are mutually exclusive,
and that the packing is hyperbolic if each degree is at least seven.
The uniform boundedness assumption was later removed by He and Schramm \cite{O3}, and they proved 
in general that the {\it type} of a packing is unique (that is, there is no infinite graph that 
packs both the disc and the plane). In the impressive paper \cite{O6},
they characterize the type in terms of the discrete extremal length, and use it to show that
the packing is parabolic if simple random walk on $G$ is recurrent. They conclude (Theorem 10.1)
that a packing is parabolic if at most finitely many vertices have degree greater than 6 
(notice that every vertex of the hexagonal packing has degree 6). This paper contains their earlier
result \cite{O6b} that a packing is hyperbolic if the lower average degree is greater than 6. By definition, the 
lower average degree is
$${\rm lav}(G)= \sup_{W_0} \inf_{W\supset W_0}  \frac1{|W|}\sum_{v\in W} {\rm deg}(v).
$$
In the case that the degrees of the vertices are uniformly bounded, they also show
that transience implies hyperbolicity. Jointly with Itai Benjamini,
this line of investigation was carried further in \cite{O12} and \cite{O13}, by  applying 
circle- and square packings to constructions of harmonic functions on graphs. 
Another nice application of circle packings is the recurrence of (weak) limits of random planar graphs
with bounded degree, \cite{O13}.

\bigskip
I have always admired Oded's ability to find a good modification of a difficult problem
that turns it into a tractable problem while keeping its essential features. 
One of the many examples is his work on discrete analytic function \cite{O14}.
Since circle packings can be viewed as discrete analogs of conformal maps, it is  natural
to ask for the analogs of analytic functions, thus giving up injectivity (disjointness of the discs).
See \cite{Ste} for the state of the art and beautiful illustrations.
Peter Doyle described collections of discs that are tangent according to the hexagonal pattern
that are analogs of the exponential function. He conjectured that these would be the only ``entire''
circle packing immersions. While Oded was not able to resolve this conjecture, he did find that 
collections of overlapping discs based on the square grid seem better suited for the problem,
and constructed the analog of the error function $\int e^{-z^2}dz$ in this setting. Along the way, he 
introduced M\"obius invariants that are discrete analogs of the Schwarzian derivative and became
instrumental in his later work \cite{O16}.
%
%
\begin{figure}[!htbp]
\begin{center}
\subfigure[]{
	\includegraphics*[viewport=1 50 600 500, width=60mm]{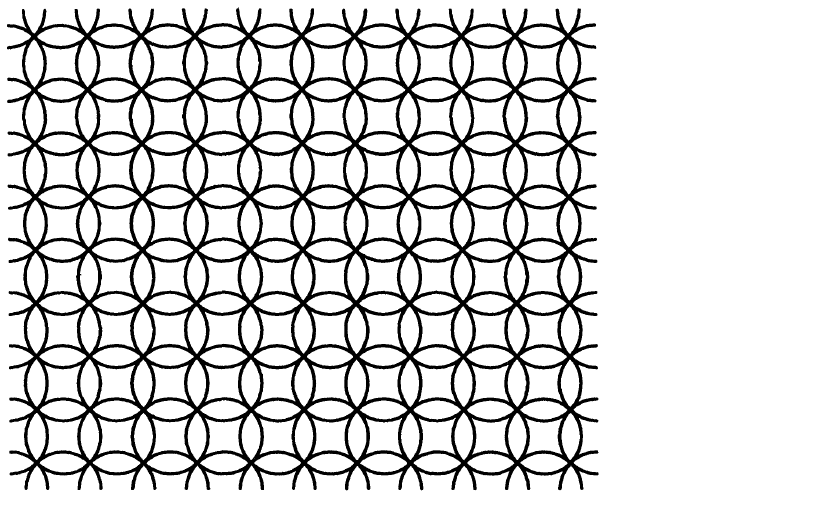}
	
}
\subfigure[]{
	\includegraphics*[viewport=10 10 1200 700, width=85mm]{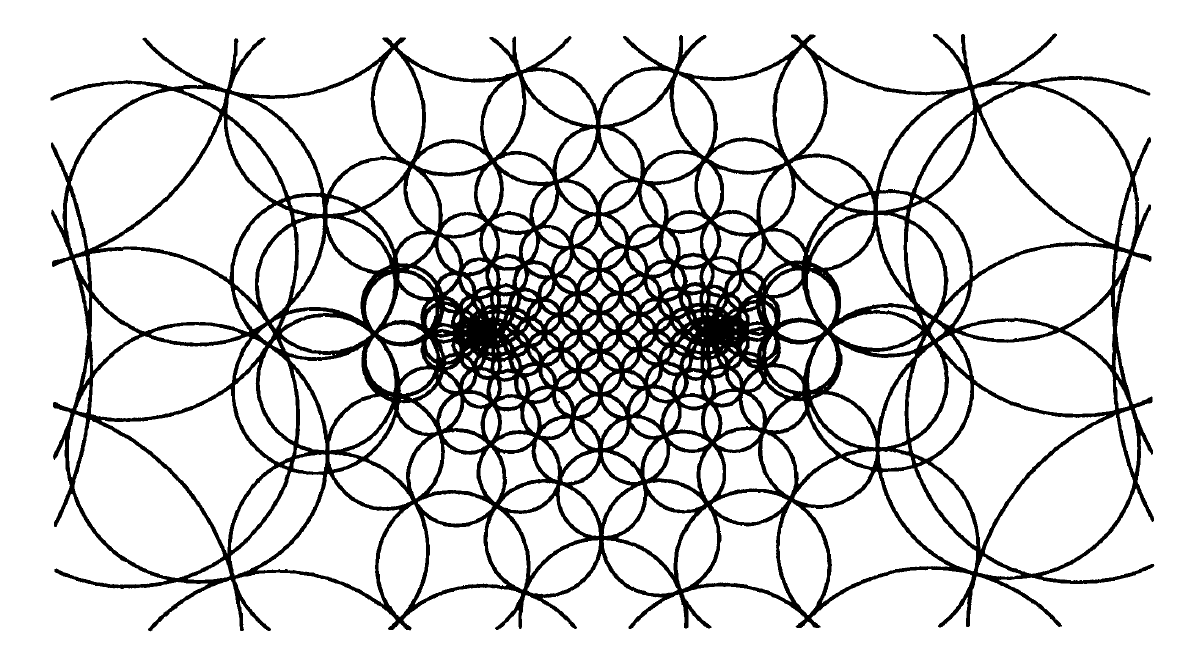}
}
\caption{The square grid and the $\sqrt{i}$SG erf pattern, from \cite{O14}}
\label{sqgrid}
\end{center}
\end{figure}
\newpage
\section{The Schramm-Loewner Evolution}

There are several excellent lecture notes, overview articles,
 and a textbook on SLE \cite{La3}, mostly by and
for probabilists  or theoretical physicists, see 
\cite{BB3}, \cite{Car2}, \cite{Dup}, \cite{GK}, \cite{KN},  \cite{S16}, \cite{W1}, \cite{W2}
and the references therein.
It is not my intention to provide another streamlined introduction to the area. 
Instead, I would like to give a somewhat historic account with an
emphasis on Oded's contributions, highlighting some of the mathematical
challenges he faced.

\subsection{Pre-history}
It is perhaps appropriate to very briefly describe the state of 
knowledge related to conformally invariant scaling limits
prior to Oded's discovery of SLE, and to describe some of the
results that were instrumental in his work. Oded's own historical narrative
is Section 1.2 in \cite{S16}.

\bigskip
Two-dimensional lattice processes such as the Self-Avoiding
Walk (SAW), the Ising model, percolation, and diffusion limited aggregation (DLA), 
to name just a few, have been intensively studied
by physicists and by probabilists for a long time. See Figure \ref{slepics}
for some pictures, and the aforementioned articles for descriptions of the models.
In the physics community, many problems such as finding the
Hausdorff dimension of scaling limits of these sets were considered
well-understood. The implicit assumption of conformal invariance 
of the scaling limit allowed the use of the powerful machinery of conformal field theory
and led to results such as Cardy's formula for the crossing probability of critical percolation
\cite{Car1}.
On the mathematical side, progress was much slower, one of the hurdles being
that in most cases the existence of a scaling limit was unknown. Even finding suitable
definitions of the concept of scaling limit was a nontrivial task. 

\bigskip
In the late 1980's, Christian Pommerenke told me how compositions of (random) conformal maps
onto slitted discs could be viewed as a variant of the Witten-Sander model for DLA \cite{WS}.
At the same time, Richard Rochberg and his son David were working on this setup. It seems that
the only trace of this is a talk given by Rochberg at the March 1990 AMS Regional 
meeting in Manhattan, Kansas, titled ``Stochastic Loewner Equation''. 
Their model is similar to an approach to Laplacian growth 
proposed  by Hastings and Levitov \cite{HL}, and is quite different from what is now called 
Stochastic Loewner Evolution or Schramm-Loewner Evolution SLE.
At that time, other analysts such as Lennart Carleson, Peter Jones and Nick Makarov 
worked with similar ideas, see e.g. \cite{CM}. Oded was at best dimly aware of these activities, 
and was not really interested in stochastic processes such as DLA until much later. 

\bigskip
Greg Lawler's invention \cite{La} of the Loop Erased Random Walk (LERW) 
provided the mathematics community with a process that shared some features with the Self-Avoiding Walk, 
but at the same time was more tractable, partly due to its Markovian property. Pemantle's work \cite{Pe}
and Wilson's algorithm provided a link between Uniform Spanning Trees (UST) and LERW's.
Intensive research on the UST \cite{Ly} culminated in the paper \cite{BLPS} by Benjamini, Lyons, Peres and Schramm. 
The deep work of Rick Kenyon \cite{Ke1}  combined powerful combinatorics and discrete complex analysis
and exhibited  conformal invariance properties of the LERW. 
He was also able to determine its expected length.

\subsection{Definition of SLE}

Oded told me in 1997 about his idea to exploit conformal invariance 
in order to study the LERW. The streamlined way to present SLE 
in courses or texts, beginning with a crashcourse on the Loewner equation
followed by a crash course on stochastic calculus (or the other way round)
is, of course, not quite representative of its emergence. In a 2006 email
exchange with Yuval Peres and myself about the history of SLE, Oded wrote:

\bigskip
{\it 
Up to the time when I started thinking about SLE, I did not really know what
Loewner's equation was, or what was the idea behind it, though I
did know that it was a tool which was important for the coefficients
problem and that it involved slit mappings and a differentiation in the
space of conformal maps. I kind of rediscovered it in the context of SLE
and then made the connection.}
\subsubsection{The (radial) Loewner equation}\label{srle}
Loewner (\cite{Lo}; see also \cite{D},\cite{P} or \cite{La3}) 
introduced his differential equation as a tool in his attempt to prove
the Bieberbach conjecture $|a_n|\leq n$ 
concerning the Taylor coefficients of normalized conformal maps $f(z)=z+\sum_{n=2}^{\infty} a_n z^n$
of the unit disc. It was also instrumental in the final solution by de Branges in 1984.

\bigskip
Let $\g$ be a simple path that is contained in $\D$ except for one endpoint on $\partial\D.$
More precisely, let $\g:[0,\tau]\to\overline\D$ be continuous and injective with $\g(0)\in\partial\D$ 
and $\g(\tau)=0$, such that $\g(0,\tau]\subset\D.$
Denote $G_t=\D\setminus\g[0,t]$ so that $G_0=\D.$ Then,
for each $0\leq t< \tau$, there is a unique
conformal map $g_t: G_t\to\D$ that is normalized by $g_t(0)=0$ and $g_t'(0)>0.$
By Schwarz' Lemma, $g_t'(0)$ strictly increases, $g_0'(0)=1,$ and it is not hard to see
that $g_t'(0)\to\infty$ as $t\to \tau.$ Hence we can reparametrize $\gamma$ so that $\tau=\infty$
and $g_t'(0)=e^{t}$. Loewner's theorem says that 
\begin{equation}\label{leg}
\frac{\partial}{\partial t} g_t(z) = g_t(z) \frac{\z_t+g_t(z)}{\z_t-g_t(z)}
\end{equation}
for all $t\geq0$ and all $z\in G_t$, where the ``driving term''
$$\z_t=g_t(\g(t))\in\partial\D$$ 
is continuous (a priori, $g_t$ is only defined in $G_t$, but it can be shown that 
$g_t$ extends to $\g(t)$).

%
%
%
\begin{figure}[!htbp]
\centering
\includegraphics*[viewport=100 220 600 720, height=80mm]{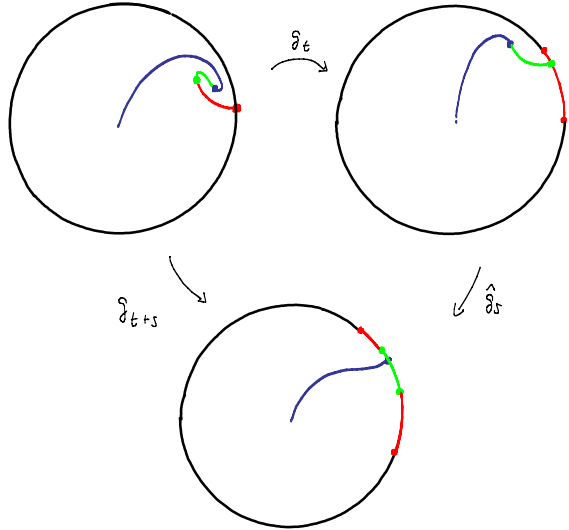}
\caption{\label{loewnerfig} Conformal maps from slit discs onto discs.}
\end{figure}

\bigskip
A simple but crucial observation is that the driving term $\z^T$ of the curve $\g^T =g_T(\gamma)$
(more precisely, the parametrized curve $\gamma^T(s) := g_T(\gamma(T+s))$)
is given by $\z^T_s = \z_{T+s}.$ Thus ``conformally pulling down'' a portion of $\g$ corresponds
to shifting the driving term. Intuitively, one can think of the Loewner equation as describing 
a conformal map to a slitted disc  as a composition of conformal maps onto infinitesimally
slitted discs with slit at $\zeta_t,$ plus the statement that the conformal map onto such a disc
is $z\mapsto z+  z \frac{\z_t+z}{\z_t-z} \Delta t$ up to first order in $\Delta t$.

\bigskip
Thus the Loewner equation associates with each simple curve $\g\subset\D$
a continuous function  $\z_t$ with values in $\partial\D.$
Conversely, it is not hard to show that the solution $g_t(z)$ to the initial value problem 
\eqref{leg}, $g_0(z)=z,$
forms a family of conformal maps of simply connected domains $G_t$ onto $\D.$ In fact, $G_t$ is 
the set of those points $\z\in\D$ for which the solution is well-defined on the interval $[0,t]$.
It easily follows that $G_t$ increases in $t,$ and that $z\in G_t$ unless $g_{s}(z)=\zeta_s$ for some 
$s\leq t.$
The complement 
$$K_t=\D\setminus G_t$$
is called the {\it hull} of $\z.$ In our original setup of a slit disc, we simply recover the curve, $K_t=\g[0,t].$
It has been known since Kufarev \cite{Ku} that smooth functions $\z$ generate smooth curves $\g,$
but that there also exist continuous functions $\z$ for which the associated hull is not a simple arc in $\D$. Kufarev's example simply is the computation that a circular chord $\g$ of the unit circle has continuous driving term.
In fact, it can be topologically wild (not locally connected), see \cite{MR}. 
My own interest in the Loewner equation originated when Oded asked me which driving terms generate curves.

\subsubsection{The scaling limit of LERW}

The LERW is obtained from simple random walk by erasing loops chronologically.
%
%
%
\begin{figure}[!htbp]
\centering
\includegraphics*[viewport=300 570 450 650, height=60mm]{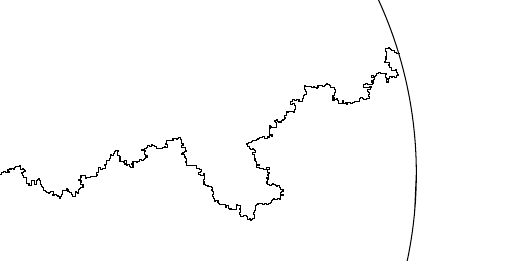}
\caption{\label{lerwfig} A loop erased random walk in the disc, from \cite{S1}.}
\end{figure}
The main result of Oded's celebrated paper \cite{S1} was a conditional 
theorem: Assuming the existence and conformal invariance of the scaling
limit of LERW, he showed that the Loewner driving term of the
resulting (random) limiting curve is a Brownian motion on the unit circle, $\z_t=e^{i B_{2t}}.$
To make this rigorous, he first gave the following definition of the notion of scaling limit:
Let $D\subsetneq\C$ be a domain, fix $a\in D$, and for $\delta>0$, consider the LERW on the graph
$\delta\Z^2\cap D$, started at a point closest to $a$ and stopped when reaching $\partial D.$
Viewing the path of the LERW as a random subset of the sphere $S^2=\C\cup\{\infty\}$, 
its distribution is a discrete measure $\mu_{\delta}$ on the space of compact subsets of $S^2$.
Equipped with the Hausdorff distance, the space of compact subsets of $S^2$ is a compact
metric space, and so is the space of its Borel measures.
The existence of  subsequential weak limits $\mu=\lim_j \mu_{\delta_j}$ follows
at once. If the limit measure 
$\mu=\lim_{\delta\to0} \mu_{\delta}$ 
exists, it is called 
the scaling limit of LERW from $a$ to $\partial D$.

\begin{thm}[\cite{S1}, Theorem 1.1]\label{curves} 
If each connected component of $\partial D$ has positive diameter, then every
subsequential scaling limit measure $\mu$ of the LERW from $a$ to $\partial D$
is supported on simple paths.
\end{thm}

In other words, the measure of the set of non-simple curves is zero.
This theorem is interesting in its own right. It has been known previously that, loosely speaking
and under mild assumptions, random curves have uniform continuity properties that imply
their (subsequential) scaling limits to be supported on continuous curves \cite{AB}.
However, the fact that the loop erased paths are {\it simple} curves does not directly imply 
that the limiting objects have no loops. Indeed, the
limits of other discrete random simple curves such as the critical percolation interface
or the uniform spanning tree Peano path are not simple. 
The proof uses estimates for the probability
distribution of ``bottlenecks'', based on harmonic measure estimates and Wilson's algorithm, 
and a topological characterization of simple curves.

\bigskip
Next, Oded formulated the 
conjecture of existence and conformal invariance of the scaling limit as follows.

\bigskip
\noindent
{\it {\bf Conjecture} (\cite{S1},1.2)\label{conjecture}
Let $D\subsetneqq\C$ be a simply connected domain in $\C$, and let
$a\in D$. 
Then the scaling limit of LERW from $a$ to $\partial D$ exists.
Moreover, suppose that
$f:D\to D'$ is a conformal homeomorphism onto a domain $D'\subset\C$.
Then $f_*\mu_{a,D}=\mu_{f(a),D'}$, where $\mu_{a,D}$ 
is the scaling limit measure of LERW from $a$ to $\partial D$,
and $\mu_{f(a),D'}$ is the scaling limit measure of LERW from $f(a)$
to $\partial D'$.
}
\bigskip

The most important and exciting result of \cite{S1} was the insight 
that this conjecture implied
an explicit construction of the limit in terms of the Loewner equation.
By Theorem \ref{curves}, the conjectural scaling limit $\mu$ induces
a measure on the space of continuous real-valued functions $\zh_t$
via the correspondence 
$\g \mapsto \z=e^{i \zh}$ of the Loewner equation. Oded showed that
the law of $\zh$ is that of a time-changed Brownian motion, $B_{2t}$:

\begin{thm}[\cite{S1}, Theorem 1.3]\label{sle} 
Assuming the above conjecture, the scaling limit $\mu$ is 
equal to the law of the hulls $K$
associated with the driving term $\z=e^{i B_{2t}}$, where $B_t,t\geq0$ is a 
Brownian motion started at a uniform random point in $[0,2\pi)$.
\end{thm}
In his characteristic way, Oded pointed out the simple idea behind the theorem.
From his paper:

\bigskip

{\it At the heart of the proof of Theorem \ref{sle} lies the following simple combinatorial fact
about LERW. Conditioned on a subarc $\b'$ of the LERW $\b$ from 0 to $\partial D$, which extends
from some point $q \in\b$ to $\partial D$, the distribution of $\b\setminus\b'$ 
is the same as that of LERW from
0 to $\partial(D - \b')$, conditioned to hit $q$. When we take the scaling limit of
this property, and apply the conformal map from $D - \b'$ to $\D$, this translates into the
Markov property and stationarity of the associated L\"owner parameter $\z$.}

He also notes that ``the passage to the scaling limit is quite delicate''. The translation 
into the Markov property and stationarity is by means of the aforementioned principle that
``conformally pulling down'' a portion $\g'$ of $\g$ corresponds
to shifting the driving term. Thus $\zh$ is a continuous process with stationary and independent
increments. Now the theory of Levy processes (and the symmetry of LERW under reflection) implies
that $\zh_t$ has the law of $\sqrt{\kappa} B_t$ for some $\kappa>0$ and a standard Brownian motion $B$.
It remained to determine the constant $\kappa.$ To this end, Oded gives the following 

\bigskip\noindent
{\bf Definition.}
{\it The (radial) stochastic Loewner evolution $SLE_{\kappa}$ with parameter $\kappa>0$
is the random process of conformal maps $g_t$
generated by the Loewner equation driven by $\z_t=e^{i \sqrt{\kappa} B_t}.$
}

\bigskip\noindent
In Section 7 of \cite{S1} he actually defined SLE as a process of random paths generated by the
Loewner equation, and therefore had to restrict to those values of $\kappa$ for which the 
resulting hulls are simple curves; he conjectured that this is the interval $[0,4]$.

Next, he analyzed the winding number of the SLE-path around 0, and of the LERW:
If $\theta_{\kappa}(t)=\arg \g(t)$, $t\geq0$, denotes the continuously defined argument 
along the curve, then he computed the variance
$$E[\theta_{\kappa}(t)^2] = (\kappa+o(1)) \log t.$$
It follows that the winding number of the portion of the SLE path until its first hitting of the
circle of radius $\eps$ centered at $0$ has variance $(\kappa+o(1)) \log(1/\eps).$
On the other hand, Kenyon's work \cite{Ke2} implies that the variance of the winding number of
LERW in $\D\cap \eps\Z^2$ is $(2+o(1)) \log(1/\eps),$ and after some work the conclusion
$\kappa=2$ follows.
\bigskip

Oded went on to compute what he called the critical value for SLE: He proved that for
$\kappa>4,$ almost surely SLE will {\it not} generate simple paths, and conjectured that
it will for $\kappa\in[0,4]$. See the discussion in Section \ref{path} for the simple
proof  using Ito's formula (in the chordal case). However, Oded was a self-taught 
newcomer to stochastic calculus, discovered some of the basics himself, and commented on his proof
in an email in January 1999:

\bigskip
{\it This must all be quite standard, to people with the right 
background.  But not for me.}
\bigskip
\subsubsection{The chordal Loewner equation, percolation, and the UST}
The classical (radial) Loewner equation is well-suited for curves that
join an {\it interior point} to a {\it boundary point}, such as curves generated by the LERW. 
Other processes generate curves joining {\it two boundary points}. 
Oded realized how important the ``correct'' 
normalizations are in dealing with conformal maps and in particular the 
Loewner equation, and found the appropriate variant of the Loewner equation (it turned out later
that this version has been described earlier, beginning with N.V. Popova \cite{Pop1},\cite{Pop2};
I would like to thank Alexander Vasiliev for this reference).
He describes this in another email in January 1999: 

%
%
%
%
\begin{figure}[!htbp]
\centering
\includegraphics*[viewport=100 530 400 750, height=50mm]{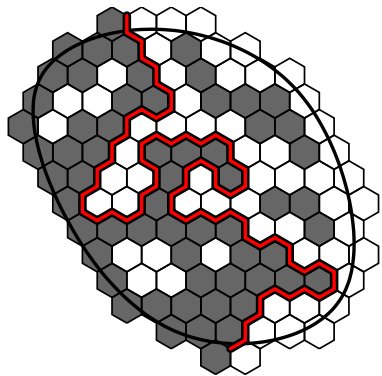}
\caption{\label{percdomainfig} Percolation in a domain, from \cite{S16}.}
\end{figure}

\bigskip
{\it 
I have a mathematical querry.
Before the question itself, here's the motivation.  For
the LERW scaling limit, the natural object is a probability measure
on the set of paths from a point in the domain to the boundary.  In
other settings, the natural object is a probability measure joining
two points on the boundary of a domain.  Consider, for example,
percolation in the unit disk.  Let $(\gamma,\beta)$ be a partition of
the boundary of the disk to two arcs, disjoint except for the endpoints.
Let $K$ be the union of all percolation clusters inside the disk that are
connected to $\gamma$.  Then the outer boundary of $K$ is a path, $\alpha$,
joining the two endpoints of $\gamma$.  The scaling limit of $\alpha$ is
conjecturally conformally invariant (but not a simple path).  Assuming
conformal invariance, I'm optimistic that the scaling limit can be
represented by a Loewner-like Brownian evolution.  The first step for this 
seems to be the following variation on Loewner's theorem:  

\bigskip\noindent
{\bf Thm:} Let $\alpha:[0,\infty)\to\C$ be a continuous simple path
such that $\alpha(0)=0, \lim_{t\to\infty} \alpha(t)=\infty$
and $\Im(\alpha(t))>0$ when $t>0.$  Let $f_t$ be the
conformal map from the upper half plane to the upper half plane
minus $\{\alpha(s): 0\leq s\leq t\}$, which is normalized by
$f_t(z) = z + O(1/z)$ near infinity.  By change of parameterization
of $\alpha$ (and perhaps changing its interval of definition),
we may assume that $f_t(z) = z + t/z + O(1)/z^2$ near infinity.
Set $g(w,t)=f_t^{-1}(w)$.  Then
 $\{\partial g/\partial t\} = 1/(g-k(t))$,
where $k(t)= g(\alpha(t),t).$

\bigskip
In the situation of percolation and related conformal invariance
models, one should expect $k(t)=c BM(t)$, where BM is on the real line.
Have you seen this theorem?  The proof should not be difficult.
It can either be derived from Loewner's theorem, or by adapting
the proof.}
\bigskip

\begin{figure}[!ht]
\begin{center}
\subfigure[]{
	\label{figa}
	\includegraphics*[viewport=10 350 600 650, height=50mm]{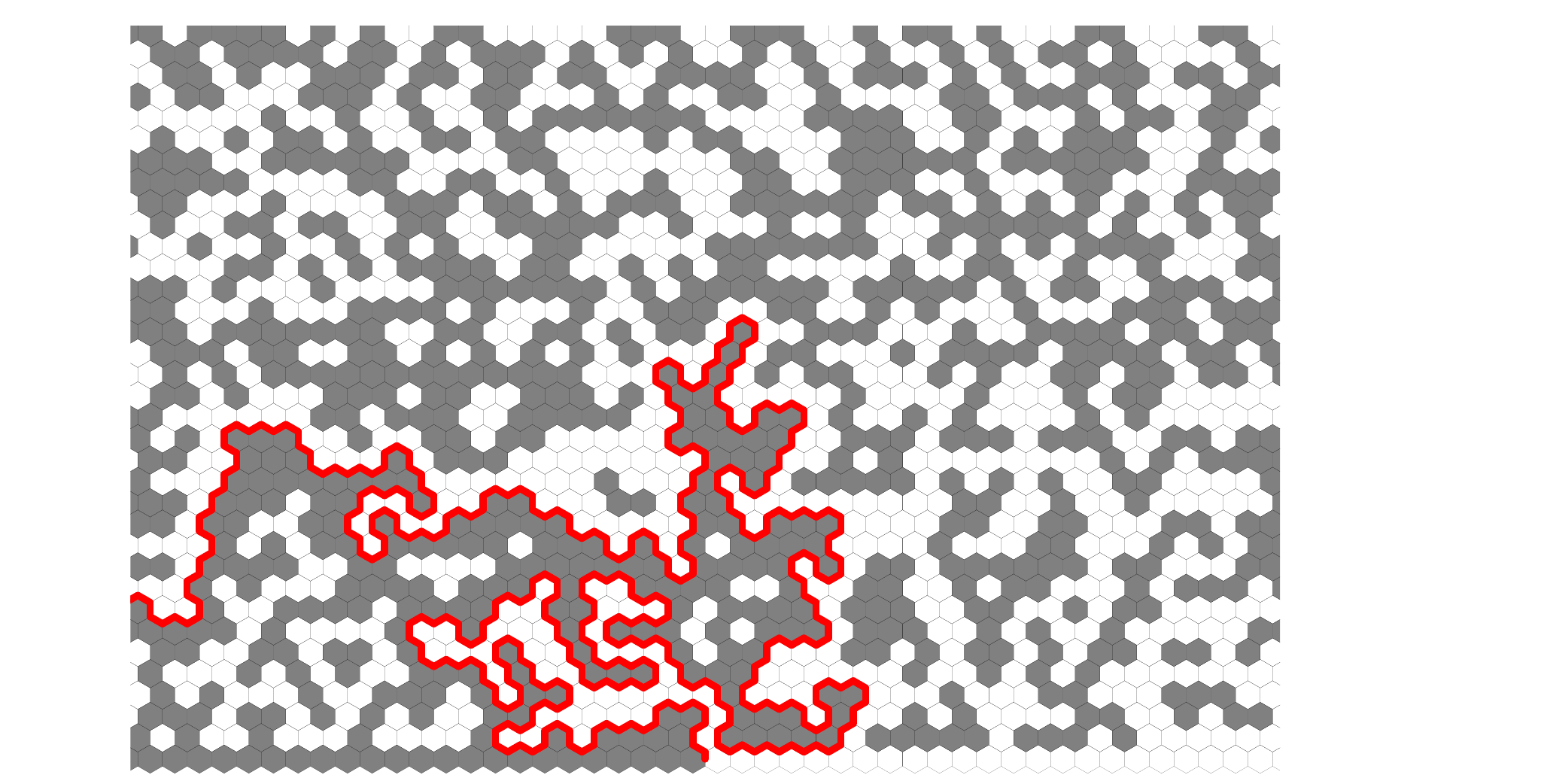}
}
\subfigure[]{
	\label{figb}
	\includegraphics*[viewport=200 350 500 650, height=60mm]{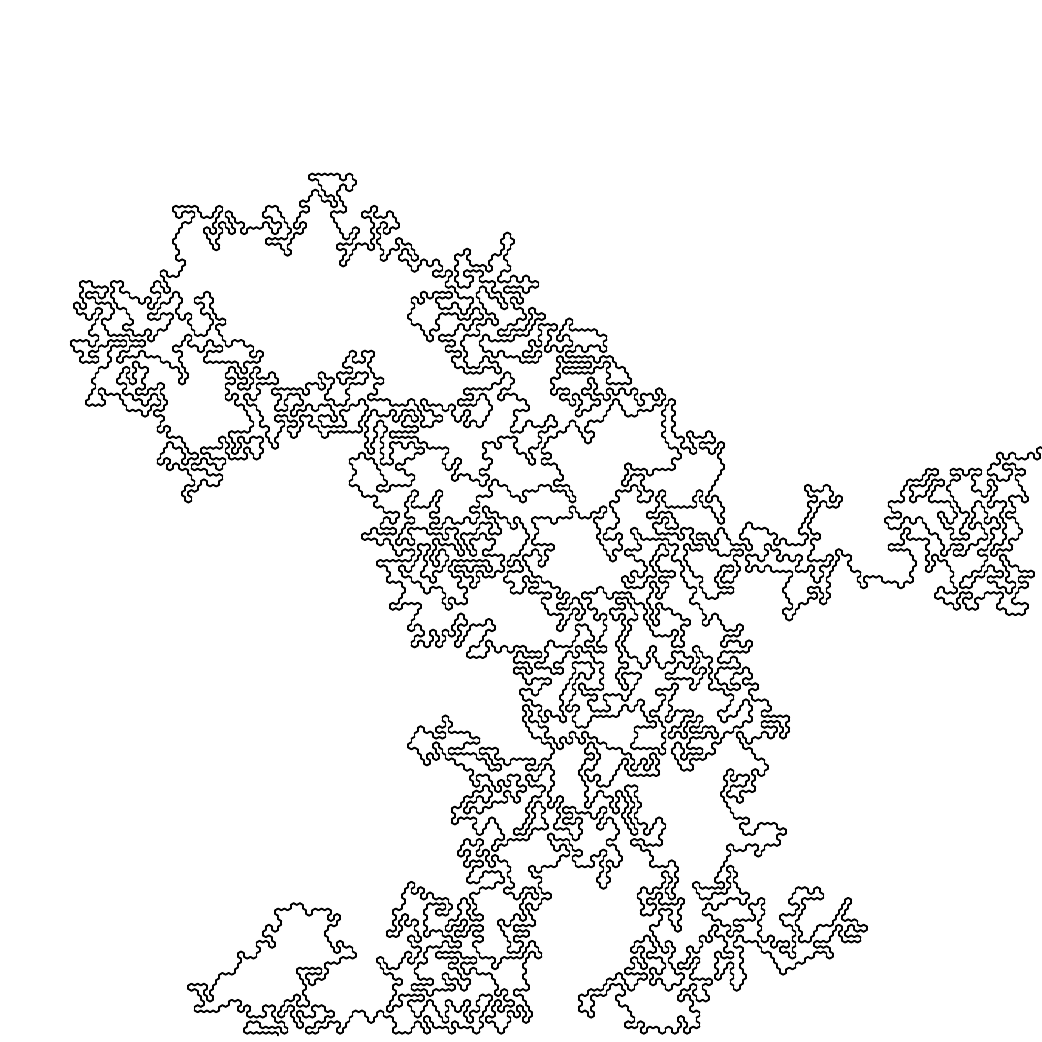}
}
\caption{The percolation exploration path}
\label{figpercolation}
\end{center}
\end{figure}

In order to coincide with the normalization of the radial Loewner equation, 
he later slightly adjusted the parametrization of the path $\alpha$  so that
$$f_t(z) = z + 2t/z + O(1)/z^2$$
and
\begin{equation}\label{cle}
\frac{\partial}{\partial t} g_t(z) = \frac{2}{g_t(z)-W_t}.
\end{equation}
With this normalization and assuming conformal invariance of the percolation scaling limit, 
he showed that
the (non-simple) limit curves would satisfy the chordal Loewner equation \eqref{cle}
with $W_t= \sqrt{6}B_t$. The value 6 can be found as the only $\kappa$ such that the random sets
generated by $\sqrt{\kappa}B_t$ satisfy Cardy's formula, or the locality property discussed below. 
This led to the definition of chordal $SLE_{\kappa}$ as the 
random process of conformal maps 
$$g_t:\H\setminus K_t\to\H$$ 
generated by the Loewner equation \eqref{cle} with driving function $W(t)=\sqrt{\kappa} B_t$, where $B$ is a standard Brownian motion.

The {\it hull} $K_t$ is the set of those points $z$ for which $g_s(z)=W_s$ for some $s\leq t$
so that \eqref{cle} becomes undefined. Since $g_t$ is determined by $K_t,$ one has the equivalent

\bigskip\noindent
{\bf Definition:} Chordal $SLE_{\kappa}$ is the process of random hulls $(K_t, t\geq0)$
generated by the Loewner equation \eqref{cle} with $W_t =  \sqrt{\kappa}B_t$.

\bigskip
In the same paper, Oded also defines and analyzes subsequential scaling limits of the uniform
spanning tree. He ends the paper by speculating (that is, stating without giving detailed proofs)
about the Loewner driving term of the UST Peano curve in the upper half plane $\H$. 
He finds that, again assuming existence and conformal invariance of the limit, 
this random space filling curve is $SLE_8$. At the end of the introduction,
he summarizes the findings of the paper as follows:

\bigskip
{\it The emerging picture is that different values of $\kappa$ in
the differential equation \eqref{leg} or \eqref{cle}
produce paths which are scaling limits of naturally defined
processes, and that these paths can be space-filling,
or simple paths, or neither, depending on the parameter $\kappa$.}

%
%
%
%
\begin{figure}[!ht]
\centering
\includegraphics*[viewport=90 100 140 120, height=60mm]{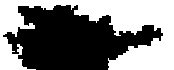}
\caption{\label{pspeanofig} The UST Peano curve, from \cite{S13}.}
\end{figure}

\subsection{Properties and applications of SLE}

Two exciting developments took place shortly after the introduction of SLE,
namely the very productive collaboration of Greg Lawler, Oded Schramm and 
Wendelin Werner, and the surprising proof of existence and
conformal invariance of the percolation scaling limit by Stas Smirnov. I will begin 
describing the former, and defer the latter to Section \ref{limit}.

\subsubsection{Locality}\label{slocality}

%
%

In the important papers \cite{LW1} and \cite{LW2}, 
Lawler and Werner discovered that Brownian excursions have a certain restriction property (explained
below), and that intersection exponents of conformally invariant processes with this 
property are closely related to those of Brownian motion. What was missing was a way to
compute exponents of {\it some} conformally invariant process.
Lawler, Schramm and Werner discovered that SLE provided such a
process to which the universality arguments of Lawler
and Werner \cite{LW2} applied. The following email from Oded describes the
crucial property: 

\bigskip
{\it
I don't remember if I've mentioned to you the restriction property for
SLE(6) that Greg, Wendelin and I have proved.  It says that (up to
time parameterization) the law
of SLE(6) in an arbitrary domain D starting from a point p on the boundary
and stopped when it exits a small ball B around p, does not depend
on the shape of the domain outside B (provided that D-B is
connected, say).  Thus, SLE(6) is purely a local process, like BM.
This is not true for $\kappa\neq6$.  For example, SLE(6) in the disk
(with Loewner's original equation) is the same as SLE(6) in the
half plane, with my variation on Loewner's equation.  This
seems to say that SLE(6) is a very special process.}

\bigskip
The above property is trivial for the discrete critical percolation exploration path, since the path
can be grown ``dynamically'' by deciding the color of a hexagon only when the path meets
it and needs to decide whether to turn ``right'' or ``left''. 
Hence, in light of the conjectured scaling limit, locality for 
chordal $SLE_6$ was not unexpected. 
The coincidence of chordal and radial 
$SLE_6$, discussed below, was more surprising.

Here is a precise statement of the locality property.
Let $D=\H\setminus A$ be a simply connected subdomain of $\H$ such that $A$ is bounded and also 
bounded away from $0.$ Denote $g_A$ the conformal map from 
$D$ onto $\H$ with
the hydrodynamic normalization ($g_A(z)-z\to 0$ as $z\to\infty$), and set $\Phi_A=g_A-g_A(0).$
Then SLE in $D$ from $0$ to $\infty$ is defined as the preimage of SLE in $\H$ from $0$ to $\infty$
under $\Phi_A$. The following expresses the fact that SLE in $D$ is a time-change
of SLE in $\H$, up to the time that the process hits $A$. Let
$T=\inf\{t: K_t\cap A\neq\emptyset\}$ and $\tilde T=\inf\{t: K_t\cap \Phi_A(\partial A)\neq\emptyset\}$.

\begin{thm}[\cite{S4}, Theorem 2.2]\label{locality} For $\kappa=6,$ the processes $(\Phi_A(K_t),t<T)$
and $(K_t, t<\tilde T)$ have the same law, up to  re-parametrization of time.
\end{thm}
The equivalence of chordal and radial $SLE_6$ was established in \cite{S5}, Theorem 4.1:
Define the hulls $K_t$ of chordal $SLE_{\kappa}$  in $\D$ from $1$ to $-1$ as the image of chordal $SLE_{\kappa}$ in $\H$, 
under the conformal map $f(z)= (i-z)/(i+z)$. Thus $K_t$ are hulls growing from 1 towards -1 in $\overline\D.$
Denote $T\leq\infty$ the first time when $K_t$ contains $0.$ 
Similarly, denote $\widetilde K_t$ the chordal $SLE_{\kappa}$ hulls in $\D$, started at $1,$ and denote
$\widetilde T$ the first time when $\widetilde K_t$ contains $-1.$ Then, for $\kappa=6,$
the laws of $(K_t,t<T)$ and $(\widetilde K_t, t< \widetilde T)$
are the same, up to a random time change $s=s(t)$. The proof and Girsanov's theorem also show that for all values of $\kappa,$ the laws of
$K_t$ and $\widetilde K_{s(t)}$ are equivalent (in the sense of absolute continuity of measures), if 
$t$ and $\widetilde t$ are bounded away from $T$ and $\widetilde T.$ See Proposition 4.2 in \cite{S5}.

\bigskip
The first proof of the locality Theorem \ref{locality} was rather long and technical, based on an analysis of the Loewner driving function of a curve 
under continuous deformation of the surrounding domain. A different and simpler proof was found later 
(see Proposition 5.1 in \cite{S10}), by analyzing
$\tilde W_t := h_t(W_t),$ where $h_t := \tilde g_t \circ g_A \circ g_t^{-1}$
and $\tilde g_t= g_{g_t(D\setminus K_t)}$ is the  normalized conformal map of $g_t(D\setminus K_t)$ to $\H$.
Writing
$$\tilde  g_t(z) = z + \frac{a_t}{z} + O(\frac1{z^2}),$$ 
computation shows that
$$\partial_t \tilde g_t(z) = \frac{\partial_t a_t}{\tilde g_t(z)-\tilde W_t} = \frac{2h_t'(W_t)^2}{\tilde g_t(z)-\tilde W_t}.$$
With $W_t=\sqrt\kappa B_t$, computation using Ito's formula shows
\begin{equation}\label{loc}
d\tilde W_t = h_t'(W_t)  dW_t + \bigl((\kappa/2) - 3\bigr)  h_t''(W_t) dt.
\end{equation}
Thus $\tilde W_t$ is a local martingale if (and only if) $\kappa=6,$ 
and  a time change shows that $(\tilde g_t, t\geq0)$ is $SLE_6$.

\bigskip

\subsubsection{Intersection exponents and dimensions}

The locality of $SLE_6$ has been used to determine the so-called {\it intersection exponents} of 2-dimensional Brownian motion, 
and to compute the Hausdorff dimensions of various sets associated with its trace. 
These results established Schramm's SLE and the Lawler-Werner universality arguments
as a fundamental and powerful new tool. 

\bigskip
If $B^1_t$ and $B^2_t$ are two independent planar Brownian motions started at two
different points $B^1_0\neq B^2_0,$ it easily follows from the subadditivity of
$t\mapsto \log \P[ B^1[0,t]\cap B^2[0,t] = \emptyset]$ that there is a number $\z>0$
such that
$$\P[ B^1[0,t]\cap B^2[0,t] = \emptyset] = \left( \frac1t \right)^{\z+o(1)}.$$

Similarly, the half-plane exponent $\tilde\z$ of the event that two independent motions  
do not intersect and stay in a halfplane is given by
$$\P[ B^1[0,t]\cap B^2[0,t] = \emptyset \ {\rm and }\ B^j[0,t]\subset\H, j=1,2] = \left( \frac1t \right)^{\tilde\z+o(1)}.$$
More generally, one considers exponents $\z_p$ for the probability of the
event that $p$ independent motions are mutually disjoint,
$\z(j,k)$ for the event that two packs of Brownian motions 
$B^1\cup\cdots\cup B^j$ and $B^{j+1}\cup\cdots\cup B^{j+k}$ are disjoint,
and the corresponding half-plane exponents $\tilde\z_p$ and $\tilde \z(j,k)$. So $\z=\z(1,1)$.
Also relevant is the {\it disconnection exponent} $2\eta_j$ for the event that the union
of $j$ Brownian motions,  started at $1$, does not disconnect $0$ from $\infty$ before time $t.$

These and other intersection exponents have been studied intensively, and
values such as $\zeta=5/8$ had been obtained by Duplantier and Kwon \cite{DK} using the mathematically non-rigorous method of conformal field theory.

An extension of $\z(j,k)$ for positive  real $k>0$ was given in \cite{LW1},
and some fundamental properties (in particular the ``cascade relations'') were established. 
In the  series of papers \cite{S4},\cite{S5},\cite{S7},\cite{S9} (see \cite{S2} for a guide and sketches of proofs), 
Lawler, Schramm and Werner were able to confirm the predictions, and they proved
\begin{thm}\label{exponents} For all integers $j\geq1$ and all real numbers $k\geq0,$
$$\z(j,k) = \frac{\sqrt{24 j+1} + \sqrt{24 k+1}-2)^2 -4}{96},\quad \z_n= \frac{4n^2-1}{48},$$
$$\tilde\z(j,k) = \frac{\sqrt{24 j+1} + \sqrt{24 k+1}-1)^2 -1}{48},\quad \tilde\z_n =\frac{2n^2+n}{6},$$
and
$$\eta_k = \zeta(k,0) =  \frac{(\sqrt{24 k+1}-1)^2-4}{48}.$$
In particular, 
$$\zeta=\frac58, \tilde\z = \frac53, \eta_1=\frac14, \eta_2=\frac23.$$
\end{thm}
The proofs are technical masterpieces combining a variety of different methods. A very rough
description is as follows: First, half-plane intersection exponents of $SLE_6$ are computed, 
based on estimates for the crossing probability of (long) rectangles. This is done by
establishing a version of Cardy's formula. Then, the universality ideas of \cite{LW2} are 
employed to pass from $SLE_6$ to Brownian motion. Finally, to cover the case $k<1$, real analyticity of 
the exponent is shown by recognizing  $e^{-2\zeta(j,k)}$ as the leading eigenvalue
of an operator $T_k$ on a space of functions on pairs of paths.

\bigskip
For a fixed time t, the {\it Brownian frontier}  is the boundary of the unbounded connected component 
of the complement of $B[0,t]$, and the set of {\it cut points} is the set of those points $p$ for which
$B[0,t]\setminus\{p\}$ is disconnected. The set of {\it pioneer points} is the union of the frontiers over all $t>0$.
Mandelbrot \cite{M} observed that the Brownian frontier 
looks like a long
self-avoiding walk. Since the Hausdorff dimension of the self-avoiding walk 
was predicted by physicists to have Hausdorff dimension 4/3, he conjectured that
the Hausdorff dimension of the Brownian frontier is 4/3. 
Greg Lawler had shown in a series of papers (see \cite{La2}) how the intersection exponents
are related to the Hausdorff dimension of subsets of the Brownian path. 
He found the values $2-2\z$, $2-\eta_2$, and $2-\eta_1$ for the dimension of the 
Brownian frontier, the set of cut points, and the set of pioneer points.
This actually required his stronger estimates of the intersection probabilites up to constant factors, rather than 
up to $(1/t)^{o(1)}$. Simpler proofs of those estimates are the content of \cite{S8}.
In combination with Theorem \ref{exponents}, this proved Mandelbrot's conjecture.

\begin{thm}[\cite{S5}, \cite{S9}]\label{dimensions} The Hausdorff dimension of 
the frontier, the set of cut points, and the set of pioneer points of 2-dimensional Brownian motion
is $4/3, 3/4$ and $7/4$ almost surely.
\end{thm}
%

%
%
%
\begin{figure}[!ht]
\centering
\includegraphics*[viewport=170 230 500 550, height=60mm]{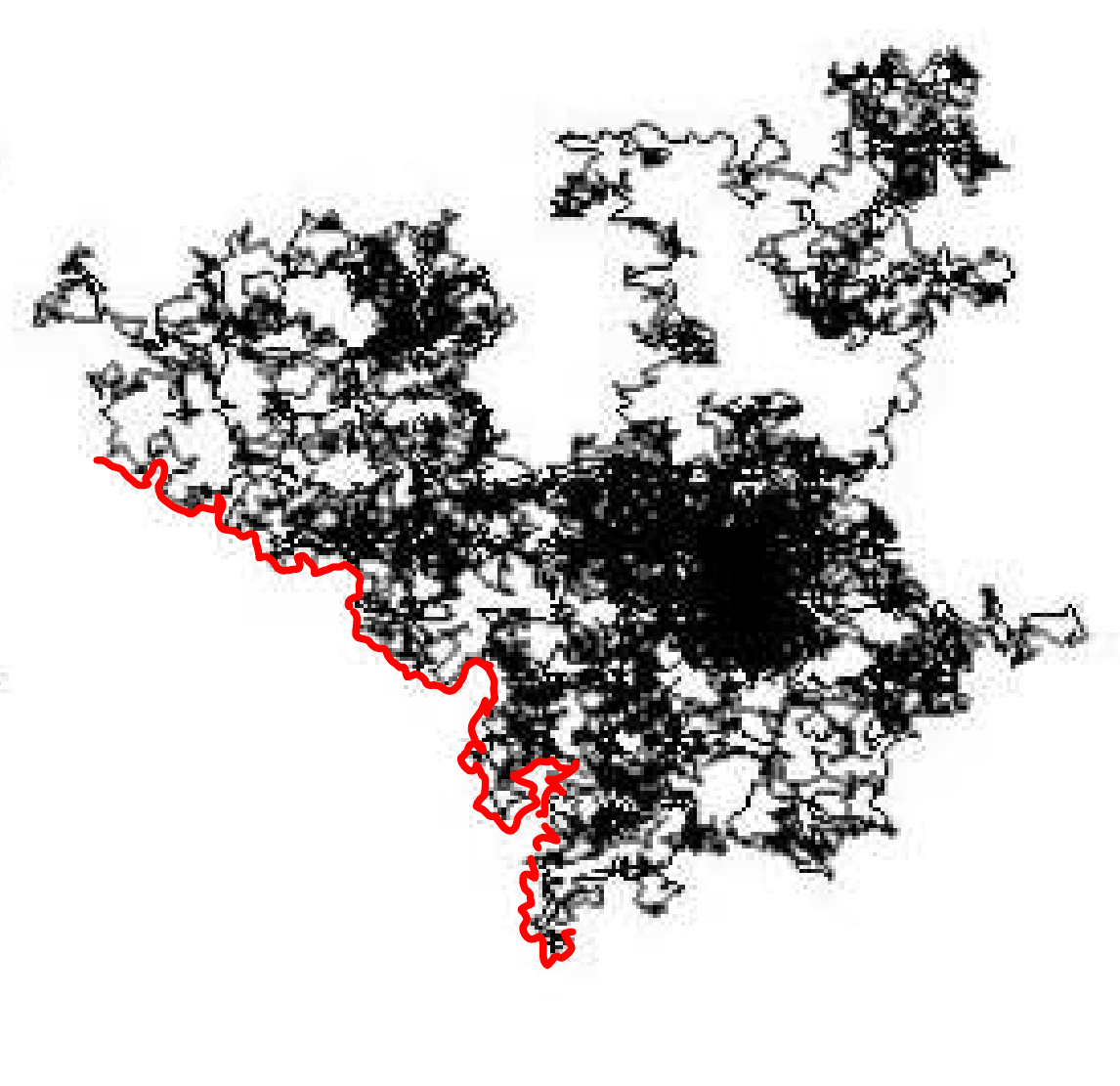}
\caption{\label{brownianfrontierfig} A Brownian path, with part of the frontier highlighted}
\end{figure}

To put this result in perspective, notice that it is rather difficult to  show even that the dimension of
the Brownian frontier is more than one \cite{BJPP}, and that the set of cutpoints is non-empty, \cite{Bu}.
It should also be mentioned that by work of Lawler, the intersection exponents for simple random walk
are the same as for Brownian motion, so that Theorem \ref{exponents} also shows, for instance, that
$$\P[S^1[0,n]\cap S^2[0,n] = \emptyset] = \left( \frac1n \right)^{\frac58+o(1)}$$
if $S_1$ and $S_2$ are two independent planar simple random walks started at different points.

Meanwhile, there is a more elegant approach to these dimension results, also due to Lawler, Schramm and Werner,
see Section \ref{srestriction}

\subsubsection{Path properties}\label{path}

I have been  fortunate  to collaborate with Oded on several 
projects over the past two decades.
Sometimes this meant just trying to catch up with 
his fast output, and watching in awe how one clever idea replaced another.
But perhaps even more impressively,
Oded had an amazing ability and willingness to listen,
and to think along. Sometimes, when I failed in an
attempt to articulate a vague idea and was about to give up a faint line of thought, 
he surprised me by completely understanding what I tried to express,
and by continuing the thought, almost like mind reading.

\bigskip
The definition of $SLE_{\kappa}$ as a family of conformal maps $g_t$ through a stochastic
differential equation does not shed much light upon the structure of the hulls 
$K_t=g_t^{-1}(\H).$ By Section \ref{slocality}, it is enough to consider chordal SLE. We say that the hull 
$(K_t)_{t\geq0}$ is {\it generated by a curve} $\gamma$ if $\gamma:[0,\infty)\to\overline\H$
is continuous and if $K_t$ is obtained from $\gamma$ by ``filling in the holes'' of 
$\gamma[0,t]$ (more precisely, $K_t$ is the complement in $\H$ of the unbounded connected 
component of $\H\setminus\gamma[0,t]$). Since continuity of the driving term $W_t$ is equivalent
to the requirement that the increments $K_{t+\eps}\setminus K_t$ have ``small diameter within $D_t$''
(more precisely, there is a set $S\subset D_t$ of small diameter that disconnects 
$K_{t+\eps}\setminus K_t$ from $\infty$ within $D_t$, see \cite{S4}, Theorem 2.6), such a curve cannot
cross itself, but it can have double points and ``bounce off'' itself (\cite{Po1}). 
There are examples of continuous $W$ for which $K_t$ is not locally connected, and such sets
cannot be generated by curves \cite{MR}. Fortunately, this does not happen for SLE:

\begin{thm}[\cite{S13}, Theorem 5.1; \cite{S11}, Theorem 4.7]\label{continuity} 
For each $\kappa>0,$ the hulls $K_t$ are generated by a curve, almost surely.
\end{thm}
It follows that, a.s., the conformal maps $f_t=g_t^{-1}$ extend continuously 
to the closed half space $\overline\H$, and $\g(t)=f_t(W_t).$ For $\kappa\neq8$, the proof
hinges on estimates for the derivative expectations $\E[|f_t'(z)|^p]$. For $\kappa=8$,
the only known proof is by exploiting the fact 
that $SLE_8$ is the scaling limit of UST, and that the UST scaling limit is a continuous curve a.s., \cite{S11}.

As Oded already noticed in \cite{S1}, the $SLE_{\kappa}$ trace has different phases,
depending on the value of $\kappa.$ 

\begin{thm}[\cite{S13}]\label{phases} 

For $\kappa\leq4,$ the SLE trace $\gamma$ is a simple curve in $\H\cup\{0\}$, almost surely. 
It ``swallows'' points (for fixed $z\in\overline\H\setminus\{0\}$, 
a.s. $z\in K_t$ for large $t$, but $z\notin\gamma[0,\infty)$) if $4<\kappa<8$, and it is space-filling ($\gamma[0,\infty)=\overline\H$)
if $\kappa\geq8$.
For all $\kappa,$ the trace is transient a.s.: $|\g(t)|\to\infty$ as $t\to\infty.$
\end{thm}
%

%
%
%
\begin{figure}[!htbp]
\centering
\includegraphics*[viewport=0 0 900 200, width=140mm]{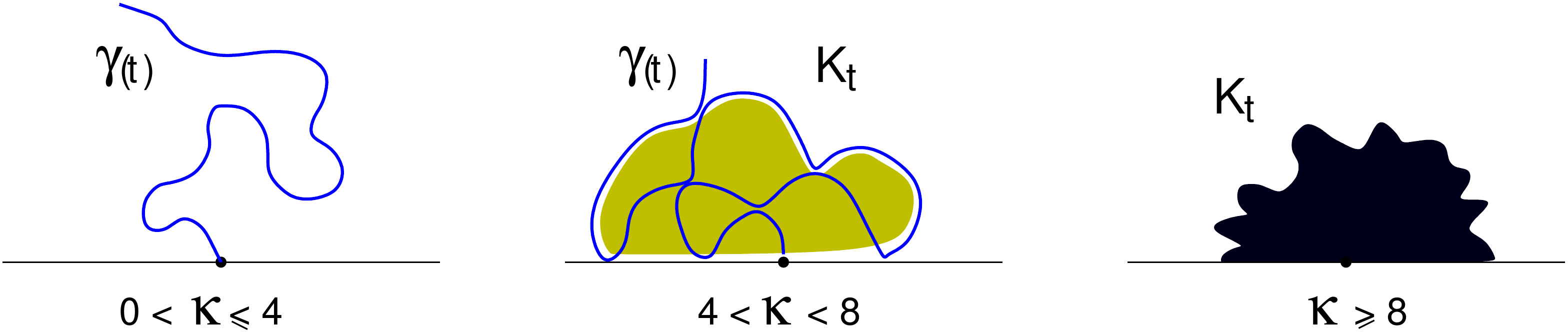}
\caption{\label{tracefig} The three phases of SLE; picture courtesy of Michel Bauer and Denis Bernard \cite{BB1}}
\end{figure}

Let us explain the phase transition at $\kappa=4$, already observed and conjectured in 
\cite{S1} (in the radial case). Let $x>0$ and denote
$$X_t=g_t(x)-W_t,$$
where $W_t=\sqrt{\kappa} B_t.$ Then
$$d X_t = \frac2{X_t}dt - \sqrt{\kappa} d B_t$$
is an Ito diffusion and can be easily analyzed using stochastic calculus. In fact,
$X_t$ is a Bessel process of dimension $1+4/\kappa$. Thus $X_t>0$ for all $t$ if and only
if $\kappa\leq 4.$ In this range, we obtain $K_t\cap\R=\{0\}$ for all $t,$ 
and it is an easy consequence of  Theorem \ref{continuity} that $\g$ is simple 
(if $r<s<t$ are such that $\g(r)=\g(t)\neq \g(s)$, then the curve
$g_s(\gamma[s,t])$ has the law of $SLE_{\kappa}$ shifted by $g_s(\gamma(s))$, but has 
two points on $\R$).

The other phase transition can be seen by examining the SLE-version of Cardy's formula:
If $X=\inf\bigl([1,\infty)\cap\gamma[0,\infty)\bigr)$ denotes the first intersection
of the SLE trace with the interval $[1,\infty),$ then a.s. $X=1$ if $\kappa\geq8$, 
whereas for $\kappa\in(4,8)$
\begin{equation}\label{cardy}
\P[X\ge s]=
\frac{
4^{(\kappa-4)/\kappa}\,\sqrt{\pi}\,
{}_2F_1(1-4/\kappa,2-8/\kappa,2-4/\kappa,1/s) \,s^{(4-\kappa)/\kappa}
}{\Gamma(2-4/\kappa)\,\Gamma(4/\kappa-1/2)}
\,,
\end{equation}
where ${}_2F_1$ denotes the hypergeometric function. At the corresponding time where $\gamma(t)=X,$
the nontrivial interval $[1,X]$ gets ``swallowed'' by $K$ at once. 
The proof of Cardy's formula in \cite{S13} is similar to the more elaborate
Theorem 3.2 in \cite{S4} and based on computing exit probabilities 
of a renormalized version of $g_t,$
$$
Y_t=\frac{ g_t(1)-W_t}{g_t(s)-W_t} \in (0,1).
$$
At the exit time $T$, we have $Y_T=0$ or $1$ according to whether $X<s$ or $>s.$ Now 
Cardy's formula can be obtained using standard methods of stochastic calculus.

\bigskip
For a simply connected domain $D\neq\C$ and boundary points $p,q$,
chordal SLE from $p$ to $q$ in $D$ is defined as the image of $SLE$ in $\H$ under a conformal map
of $\H$ onto $D$ that takes $0$ and $\infty$ to $p$ and $q$.
Since the conformal map between $\H$ and $D$  generally does not extend to $\overline\H$,
the continuity of the SLE trace in $D$ does not follow from  Theorem \ref{continuity}.
However, using Theorem \ref{tdimension} below and general properties of conformal maps, 
it can be shown to still hold true, \cite{S21}.
Another natural question is whether SLE is {\it reversible}, namely if SLE in $D$ from $p$ to $q$
has the same law as SLE from $q$ to $p.$ This question was recently answered positively for $\kappa\leq4$ 
by Dapeng Zhang \cite{Z1}. It is known to be false for $\kappa\geq8$ \cite{S13}, and unknown
for $4<\kappa<8.$
\begin{thm}[\cite{Z1}]\label{reversability} 
For each $\kappa<4,$ $SLE_{\kappa}$ is reversible, and for $\kappa\geq 8$ it is not reversible.
\end{thm}

\bigskip
The aforementioned derivative expectations $\E[|f_t'(z)|^p]$  also led to 
upper bounds for the dimensions of the trace and the frontier.   
The technically more difficult lower bounds were proved by Vincent Beffara \cite{Be}
for the trace.

For $\kappa>4$, notice that the
outer boundary of $K_t$ is a simple curve joining two points on the real line.
There is a relation between $SLE_{\kappa}$ and $SLE_{\16/\kappa}$, first derived by Duplantier 
with mathematically non-rigorous methods, and recently proved in the  papers of Zhang \cite{Z2} 
and Dubedat \cite{Dub2}.
Roughly speaking, Duplantier duality says that this curve is $SLE_{\16/\kappa}$ between
the two points. A precise formulation is based on a generalization of $SLE$, the so-called
$SLE(\kappa,\rho)$ introduced in \cite{S10}. As a consequence, the dimension of the frontier can 
thus be obtained from the dimension of the dual SLE.

Based on a clever construction
of a certain martingale, in \cite{S20} Oded and Wang Zhou determined the size of the intersection of the trace with the real line. The
same result was found independently and with a different method by Alberts and Sheffield \cite{AlSh}. Summarizing:

\begin{thm}\label{tdimension} For $\kappa\leq 8,$ 
$$\rm{ dim }\ \g[0,t]=1+\frac{\kappa}{8}.$$
For $\kappa>4,$
$$\rm{ dim }\ \partial K_t = 1+\frac2{\kappa}.$$
For $4<\kappa<8,$  
$$\rm{ dim }\ \g[0,t]\cap\R = 2-\frac8{\kappa}.$$
\end{thm}
\noindent
The paper \cite{S20} also examined the question how the SLE trace tends to infinity. Oded and Zhou  showed
that for $\kappa<4$, almost surely $\g$ eventually stays above the graph of the function
$x\mapsto x(\log x)^{-\beta},$ where $\beta=1/(8/\kappa-2).$

\subsubsection{Discrete processes converging to SLE}\label{limit} 
In \cite{S4}, Lawler, Schramm and Werner wrote that 

{\it ... at present, a proof of the conjecture that $SLE_6$ is the scaling limit of critical percolation cluster boundaries seems out of reach...}

Smirnov's  proof \cite{Sm1} of this conjecture
came as a surprise. More precisely, he proved convergence of the critical site percolation
exploration path on the triangular lattice (see Figure \ref{figb})  to $SLE_6$. See also \cite{CN} and \cite{Sm2}.
This result was the first instance of a statistical physics model proved to converge to an SLE. 
The key to Smirnov's theorem is a version of Cardy's formula.
Lennart Carleson realized that Cardy's formula assumes a very simple form when
viewed in the appropriate geometry: When $\kappa=6,$ the right hand side $f(s)$ of \eqref{cardy} 
is a conformal map of the upper half plane onto an equilateral triangle $ABC$ such that $0,1$
and $\infty$ correspond to $A,B$ and $C.$ Since $SLE_6$ in  $ABC$ from $A$ to $B$ has the same law 
as the image of 
$SLE_6$ in $\H$ from $0$ to $\infty,$  the first point $X'$ of intersection with $BC$
has the law of $f(X).$ It follows that $X'$ is uniformly distributed on $BC.$
(A similar statement is true for all $4<\kappa<8$, where ``equilateral'' is replaced by ``isosceles'', and
the angle of the triangle depends on $\kappa$, \cite{Dub1}).
Smirnov proved that the law of a corresponding observable on 
the lattice converges to a harmonic function, as the lattice size tends to zero. And he was able to
identify the limit, through its boundary values.  
The proof makes use of the symmetries of the triangular
lattice, and does not work  on other lattices such as the square grid, where convergence is still unknown.

\bigskip
The next result concerning convergence to SLE was obtained by the usual suspects Lawler, Schramm and Werner \cite{S11}. 
They proved Oded's original Conjecture \ref{conjecture} about convergence of  LERW 
to $SLE_2$, and the dual result (also conjectured in \cite{S1}) that the UST converges to $SLE_8$, see Figures
\ref{lerwfig} and \ref{pspeanofig}. 

\begin{figure}[!ht]
\centering
\includegraphics*[viewport=0 0 370 200, height=60mm]{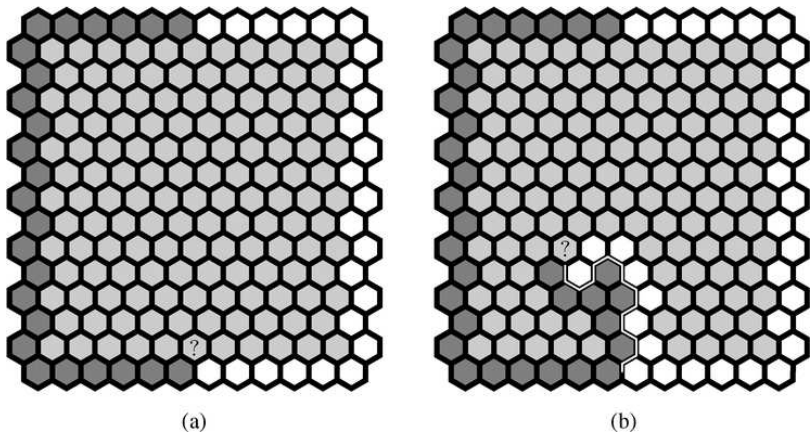}
\caption{\label{harmexplfig1} Definition of the Harmonic explorer path, from \cite{S14}.}
\end{figure}

\bigskip
The {\it harmonic explorer} is a (random)  interface defined as follows: Given a planar simply connected domain with two marked boundary points that partition the boundary into black and white hexagons, color all hexagons in the interior of 
the domain grey, see
Figure \ref{harmexplfig1} (a). 
%
%
The (growing) interface $\gamma$ starts at one of the marked boundary points and keeps the black hexagons on its left and the white hexagons 
on its right. It is (uniquely) determined (by turning left at white hegagons and right at black) until a grey hexagon is met. When it meets a grey hexagon $h$ (marked by ? in Figure \ref{harmexplfig1})
the (random) color of $h$ is determined as follows.  
A random walk on the set of hexagons is started, beginning with the hexagon $h$. The walk stops as soon as it meets a white or black hexagon, and $h$ assumes that color. Continuing in this fashion, $\gamma$ will eventually reach the other boundary point. In \cite{S14},
Oded and Scott Sheffield showed (distributional) convergence of $\gamma$ to $SLE_4$. 
%
%
%
%
\begin{figure}[!ht]
\centering
\includegraphics*[viewport=0 60 300 230, height=60mm]{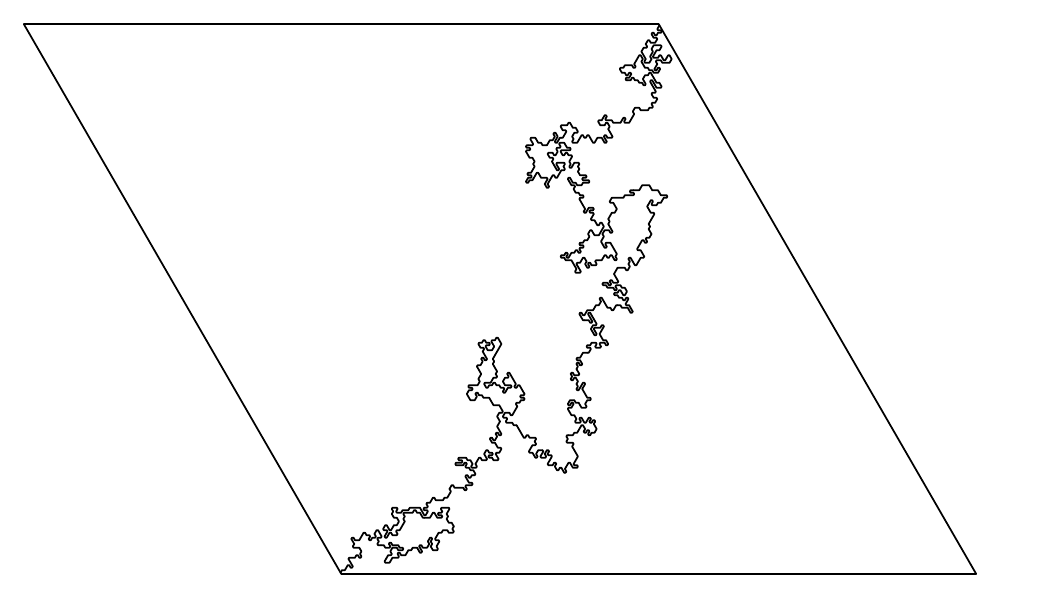}
\caption{\label{harmexplfig2} Harmonic explorer path, from \cite{S16}.}
\end{figure}
The overall strategy is again to directly analyze the Loewner driving
term of the discrete path. The crucial property of $SLE_4$ is that,
conditioned on the $SLE$ trace $\gamma[0,t]$, the probability that a point $z\in\H$ will end
up on the left of $\gamma[0,\infty)$ is a harmonic function of $z$ 
(it is equal to the argument of  $g_t(z)-W_t$, divided by $\pi$).
 
Other processes are believed to converge to $SLE_4,$ too, in particular 
Rick Kenyon's  double domino path, and the $q-$state Pott's model with $q=4.$

\bigskip
The self-avoiding walk, first proposed in 1949 as a simple model for the structure of polymers, 
has played an important role in the development of SLE, in several ways: First, Lawler's invention
of the LERW was partly motivated by the desire to create a  model that is simpler than SAW. Second, the apparent
similarity to the Brownian frontier motivated Mandelbrot's conjecture. 
Third, and most significantly, the SAW is conjectured to converge to $SLE_{8/3}.$
See \cite{S12} for precise formulations, and a proof of this conjecture assuming existence and conformal invariance
of the scaling limit, and \cite{K} for  strong numerical evidence.
However,  still very little is known rigorously about the SAW.
\begin{figure}[!ht]
\centering
\includegraphics*[viewport=320 70 670 350, height=60mm]{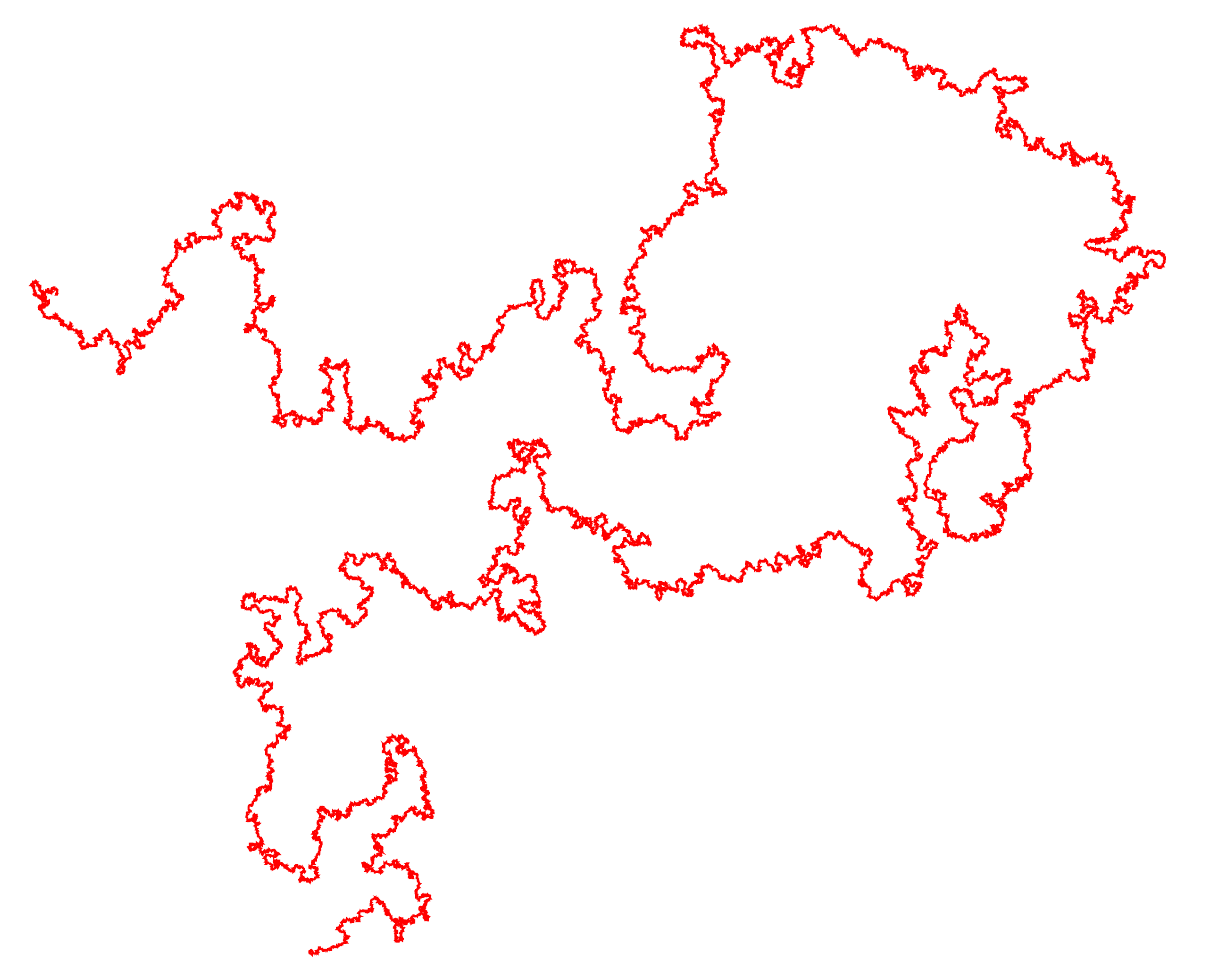}
\caption{\label{sawfig} half-plane SAW, picture courtesy of Tom Kennedy.}
\end{figure}

\bigskip 
Another famous classical  model is the Ising model for ferromagnetism. 
Stas Smirnov \cite{Sm3} has recently obtained another breakthrough concerning
convergence of lattice models to SLE. He found observables for
the Ising model at criticality and was able to prove their conformal invariance in the scaling limit. As a consequence, he obtained
$SLE_3$ in the limit. Quoting from \cite{Sm2}:

{\bf Theorem.} {\it As the lattice step goes to zero, interfaces in Ising and Ising random cluster
models on the square lattice at critical temperature converge to SLE(3) and SLE(16/3) correspondingly.}
\begin{figure}[!ht]
\centering
\includegraphics*[viewport=0 20 670 670, height=60mm]{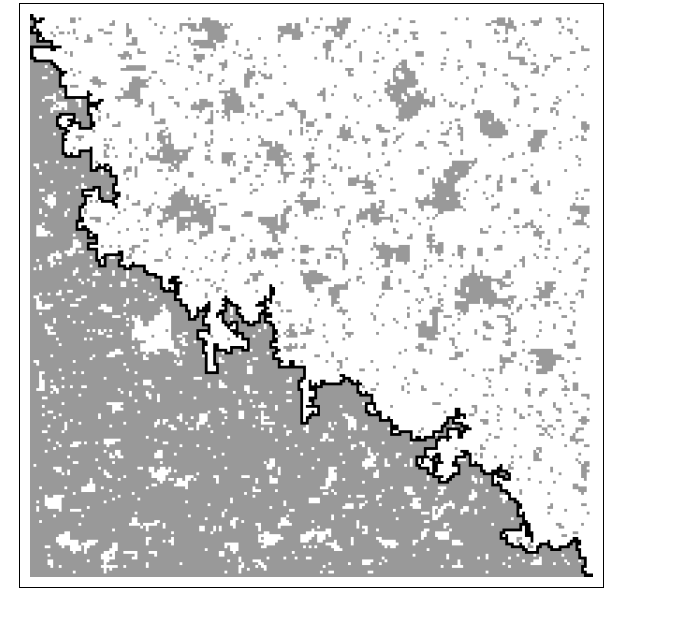}
\caption{\label{critisingfig} Critical Ising interface, picture courtesy of Stas Smirnov \cite{Sm2}.}
\end{figure}

\bigskip
\subsubsection{Restriction measures}\label{srestriction}  
The elegant and important paper \cite{S10} is a culmination of the 
universality arguments that have been initiated in \cite{LW2} 
and developed in the subsequent collaboration of Lawler, Schramm and Werner.
In the setting of random sets joining two boundary points of a simply connected
domain,  \cite{S10} gives a complete characterization of laws satisfying the
conformal restriction property, and various constructions of them.

Roughly speaking, a family of random sets $K$ joining 0 and $\infty$ in $\H$
satisfies conformal restriction, if for every reasonable
subdomain $D=\H\setminus A$ of $\H$, 
the law of $K$ conditioned on $K\subset D$ is the same as
the law of $g(K),$ where $g$ is a conformal map from $D$ to $\H$ fixing $0$
and $\infty.$ More precisely, the sets $K$ are supposed to be connected, have connected complement, and are such that 
$\H\setminus K$ has two connected components. The subdomain $D$ is reasonable
if it is simply connected and contains (relative) neighborhoods of $0$ and $\infty$. 

An equivalent definition is to consider, for each simply connected domain $D$ and
each pair of boundary points $a,b$, a law $P_{D,a,b}$ on subsets of $D$ joining $a$
and $b.$ Then the two required properties are conformal invariance, namely
$$g_* P_{D,a,b} = P_{g(D),g(a),g(b)}$$
for conformal maps $g$ of $D$,
and ``restriction'': For reasonable $D'\subset D,$ 
the law $P_{D,a,b}$ of $K$, when restricted to $K\subset D'$, equals $P_{D',a,b}$.
The remarkable main result is that there is a unique one-parameter family of such measures.
\begin{thm}\label{restriction} 
$\P=P_{\H,0,\infty}$ is a conformal restriction measure if and only if there is $\alpha>0$
such that
\begin{equation}\label{eqrestriction}
\P[K\subset D] = g_D'(0)^{\alpha}
\end{equation}
for every reasonable $D\subset\H.$
For each $\alpha\geq\frac58$ there is a conformal restriction measure $P_{\alpha}$.
Furthermore,  $\alpha_0=\frac58$ is the smallest $\alpha$ for which there is a restriction measure, $P_{5/8}$ is the only restriction measure supported on simple curves, and $P_{5/8}$ is $SLE_{8/3}.$
\end{thm}
An important observation, due to Balint Virag \cite{V},  
is that  $P_1$ is the law of Brownian excursions from $0$ to $\infty$ in $\H$
(roughly, Brownian motion started at $0$ and conditioned to ``stay in $\H$'' for all time).
An elegant application goes as follows.
If $K_1$ and $K_2$ are independent samples from $P_{\alpha_1}$ and $P_{\alpha_2}$,
then \eqref{eqrestriction} implies that
$K_1\cup K_2$ has the law of $P_{\alpha_1+\alpha_2}$ (after the ``loops'' of the union
have been filled in ). By uniqueness, it follows that the law of the union of
5 independent Brownian excursions in $\H$ (plus loops) is the same as that of 
8 copies of $SLE_{8/3}$ (with loops added). In particular, the frontiers are the same 
and thus have Hausdorff dimension $4/3$ by Theorem \ref{tdimension}.
A similar result is \cite{S10}, Theorem 9.1: The law of the hull of whole-plane $SLE_6$, 
stopped when reaching the boundary of a disc $D$, is the same as the law of a planar Brownian
motion (with the bounded complementary components added), stopped upon leaving $D.$

\bigskip
The proof that $SLE_{8/3}$ is $P_{5/8}$ is based on the following computation.
Using the same notation as \eqref{loc}, one can show

\begin{equation}
d\ h_t'(W_t) = h_t''(W_t) dW_t + \Bigl(\frac{h_t''(W_t)^2}{2 h_t'(W_t)}  
+ \bigl(\frac \kappa 2 - \frac 4 3 \bigr) h_t'''(W_t)\Bigr)dt .
\end{equation}
For $\kappa=\frac83,$ it follows that
$$ d\ h_t'(W_t)^{5/8} = \frac58 \frac{h_t''(W_t)}{h_t'(W_t) ^{3/8}} dW_t$$
so that $h_t'(W_t)^{5/8}$ is a local martingale. Writing as before 
$T=\inf\{t: K_t\cap A\neq\emptyset\}$, it is not hard to show that
$h_t'(W_t)$ tends to 0 as $t\to T$ if $K\cap A\neq\emptyset$ (the case $T<\infty$),
and $\lim_{t\to T}h_t'(W_t)\to1$ otherwise. Thus 
$$\P[K\subset D] = \P[T=\infty] = \E[h_T'(W_T)^{5/8}] = \E[h_0'(W_0)^{5/8}] = g_D'(0)^{5/8}.$$
For values greater than $5/8,$ there are several constructions of the restriction measures 
described in \cite{S10}. One is by adding ``Brownian bubbles'' to SLE-traces.

\subsubsection{Other results}

There are other versions of the Loewner equation. The ``whole plane'' equation was developed 
and used in \cite{S5} to deal with hulls $K_t$ that are growing in the plane rather than a disc
or half-plane. ``Di-polar SLE'' was introduced in \cite{BB2}, see also \cite{BBH}.
An important generalization of SLE are the $SLE(\kappa,\rho)$ and variations, 
first introduced in \cite{S10}. An elegant and unified treatment of all these 
variants is in \cite{S15}.
Also, defining SLE in multiply connected domains creates a new difficulty that is not present in the
simply connected case, since a slit multiply connected domain  is not conformally equivalent to
the unslit domain. See \cite{Z},\cite{BF1},\cite{BF2}.

\bigskip

Since SLE is amenable to computations, the convergence of discrete processes to SLE can be used to 
obtain results about the original process. In this fashion, Oded \cite{S3} obtained the 
limiting probability, as the lattice size tends to zero in critical site percolation on the triangular 
lattice in the disc $\D$, that the union of a given arc $A\subset\partial\D$ and a percolation cluster
surrounds $0$. In \cite{S6}, 
Lawler, Schramm and Werner showed that the probability  of the event $0\leftrightarrow C_R$
that the percolation cluster containing the origin reaches the circle of radius $R$ behaves like $R^{-5/48},$ 
$$P[0\leftrightarrow C_R] = R^{-5/48 + o(1)}$$
as $R\to\infty.$ See also \cite{SW} for related exponents.

\bigskip

Because of space, in this note we have ignored the mathematically nutritious ``Brownian loop soup'' \cite{LW3} and its relation to restriction measures, as well as the 
growing literature around the important Conformal Loop Ensemble $CLE_{\kappa}$ introduced
by Scott Sheffield \cite{Sh2}. See \cite{W3} and \cite{S18}.

\bigskip
There are deep and exciting connections between the Gaussian Free Field and SLE, as explored
by Oded and Scott Sheffield. The GFF has made its first appearance in this area in Rick Kenyon's 
work on the height of domino tilings \cite{Ke3}. See \cite{Sh1} for definitions and properties.
Here is a very brief description of their work.
Let $D\subset\C$ be a domain bounded by a simple closed curve that is partitioned into two arcs
by two marked boundary points.
Approximate $D$ by a portion $G=(V,E)$ of the triangular grid as before (see Figure \ref{DGFF}, where again vertices are
represented by hexagons), and denote $\partial V$ the boundary vertices. Fix a constant $\lambda$, and 
let $h=h_{\eps}$ be an instance of the Discrete Gaussian Free Field, with boundary values $\pm\lambda$
on the two boundary arcs.
This means that $h(v), v\in V\setminus \partial V$, is a $\sharp(V\setminus \partial V)$-dimensional Gaussian random variable 
whose density is proportional to $\exp(-\sum_{(u,v)\in E} (h(v)-h(u))^2/2).$ 
Extend $h$ in a piecewise linear fashion from the vertices to the triangles.
The main result of the deep and very long paper \cite{S17} is, roughly speaking, the following.
If $\lambda=3^{-1/4}\sqrt{\pi/8},$ 
then the level curve $\gamma_{\epsilon}$ of level $h=0$, joining the two marked boundary points, converges to $SLE_4$
as $\eps\to0.$ Other values of $\lambda$ lead to variants of $SLE_4.$ See also \cite{S17b} and {\cite{Dub3}.

\begin{figure}[!ht]
\centering
\includegraphics*[viewport=100 150 570 450, height=80mm]{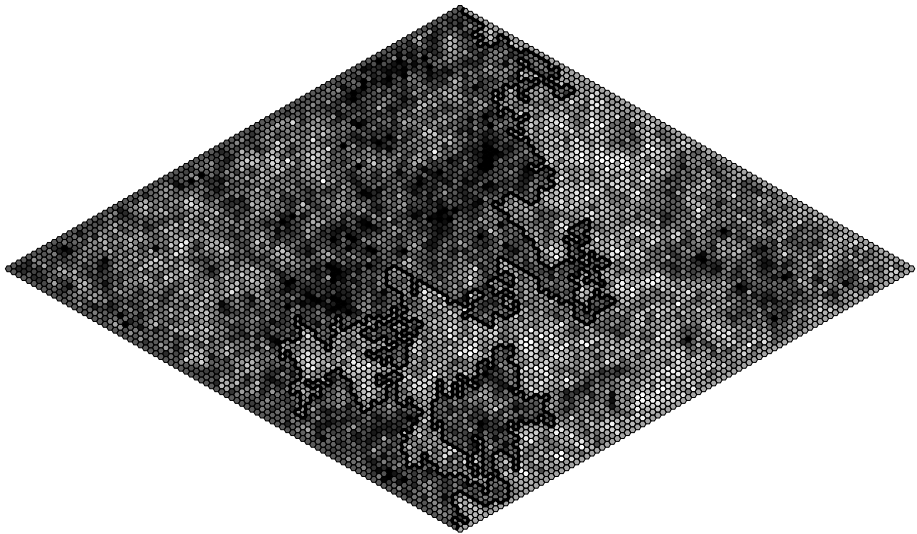}
\caption{\label{DGFF} Level set of the Discrete GFF, from \cite{S17}.}
\end{figure}

\bigskip
Finally, there are several collaborations of Oded being written at the moment. For instance,
there are deep results of Christophe Garban, Gabor Pete and Oded concerning near-critical 
percolation and its scaling limit, which is different from $SLE_6$, see \cite{GPS}.
The nice paper
\cite{ShWi} describes Oded's (unpublished) proof of Watt's formula for double crossings in critical
percolation, and provides insight into Oded's masterful use of Mathematica. 
Watt's formula was first proved rigorously by Dubedat \cite{Dub06}.

\subsection{Problems}

Many of Oded's papers contain open problems, some (such as \cite{O9}) even propose
a direction to tackle them. His ICM talk \cite{S16} contains a large number
of problems around SLE, and has provided the field with a sense of direction.
Some additional SLE-related problems are in \cite{S13}. 
Several of his problems have been solved
since their publication.

\bigskip
As already mentioned, 
the convergence of the Ising interface to $SLE_3$, Problem 2.5 in \cite{S16},
was proved by Smirnov.
The reversability of the chordal SLE path, Problem 7.3 of \cite{S16},
has been established for $\kappa\leq4$ by Zhang \cite{Z1}, and the method has even
led to a proof of a version of Duplantier's duality by Zhang and by Dubedat, \cite{Z2}, \cite{Dub2}.
Near-critical percolation and dynamical percolation  (Sections 2.6 and 5 of \cite{S16}) 
are now well-understood by \cite{GPS} and \cite{NW}.

\bigskip
Spectacular progress has been made concerning the mathematical foundations
of quantum gravity. 
Whereas the existence of an (unscaled) limit of ``random triangulations of the sphere'' 
was already established in \cite{AnS}, the important Knizhnik-Polyakov-Zamolodchikov formula
(relating exponents of statistical physics
models in ``random geometries'' to corresponding exponents in plane geometry)
seemed out of reach until recently. In Section 4 of \cite{S16}, Oded wrote:

\bigskip
{\it However, there is still no mathematical understanding of the KPZ formula. In fact, the
author's understanding of KPZ is too weak to even state a concrete problem.}
\bigskip

The Problem 4.1 in \cite{S16},  to show the existence of the (weak) Gromov-Hausdorff scaling limit of
the graph metric on random triangulations of the sphere, 
was solved in the impressive work of Le Gall \cite{LG}. Le Gall and Paulin showed \cite{LGP}
that the limiting space is a topological sphere, almost surely. Duplantier and Sheffield \cite{DS}
described a random measure (a scaling limit of the measure $e^{c h_{\eps}} dxdy$ where  $h_{\eps}$
is the Gaussian free field, averaged over circles of radius $\eps$) which exhibits a KPZ-like
relation. They conjecture a precise relation between a scaling limit of \cite{AnS} and their random
continuous space, and discuss connections to SLE. 
Following Duplantier and Sheffield, simpler random metric spaces exhibiting
KPZ were considered in \cite{O17} and \cite{RV}.

\subsection{Conclusion}
We have seen how Oded shaped the field of circle packings, and how he 
developed a deep understanding of discrete approximations to conformal maps.
His results on the Koebe conjecture are still the best to date.
We have also seen how Oded's discovery of SLE led to a powerful new tool in 
probability theory and in mathematical physics. In fact, it has changed the
way physicists and mathematicians think about critical lattice interfaces,
and has led to very fruitful interactions across disciplines.
The number of mathematicians and physicists working with SLE is increasing
fast, and the last few years have seen a number of exciting developments.
Oded has already established his place in the history of mathematics.
I have no doubt that we will see many more wonderful developments directly or indirectly 
related to Oded's work, thus keeping his spirit alive through the work of his fellow 
mathematicians, coauthors, and friends.

\begin{figure}[!htbp]
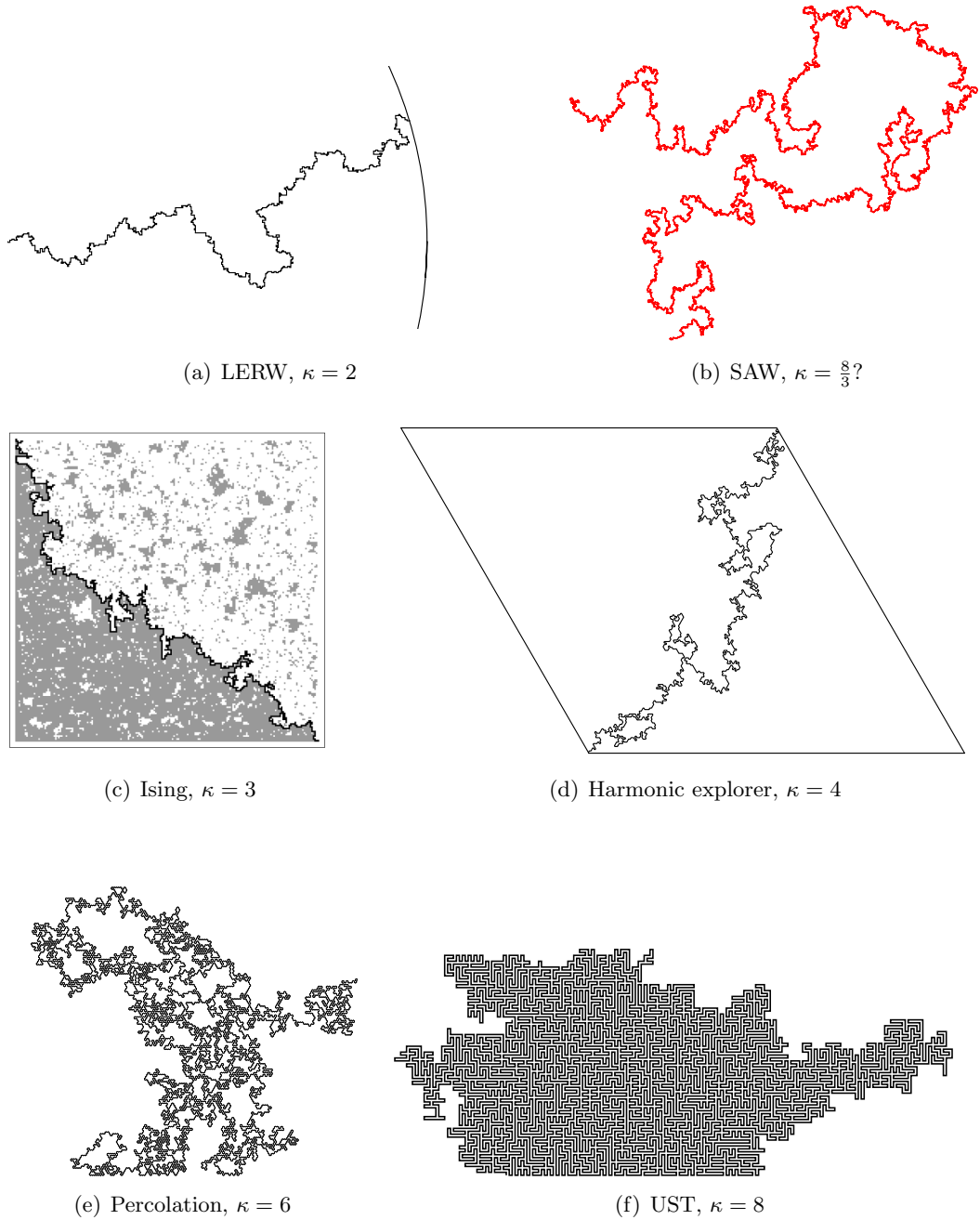

\begin{center}
\subfigure[LERW, $\kappa=2$]{
	\includegraphics*[viewport=300 570 450 650, height=40mm]{./pics/lerw.pdf}
}
\subfigure[SAW, $\kappa = \frac83$?]{
	\includegraphics*[viewport=320 70 670 350, height=50mm]{./pics/halfb.pdf}
}
\subfigure[Ising, $\kappa=3$]{
	\includegraphics*[viewport=0 20 670 670, height=50mm]{./pics/ising.png}
}
\subfigure[Harmonic explorer, $\kappa=4$]{
	\includegraphics*[viewport=0 60 300 230, height=50mm]{./pics/hepathsimulation.pdf}
}

\subfigure[Percolation, $\kappa=6$]{
	\includegraphics*[viewport=200 350 500 650, height=50mm]{./pics/Lpcurve.pdf}
}
\subfigure[UST, $\kappa=8$]{
	\includegraphics*[viewport=90 100 140 120, height=35mm]{./pics/pspeano.pdf}
}
\caption{\label{slepics} Various random curves converging to SLE's}
\end{center}
\end{figure}

\end{document}